\documentclass[11pt]{article}

\usepackage[left=2.7cm,bottom=2.7cm,right=2.7cm,top=2.7cm]{geometry}

\usepackage{color,hyperref}

\usepackage[numbers]{natbib}
\newif\ifworkinprogress
\workinprogresstrue

\ifworkinprogress
\newcommand{\RM}[1]{\textcolor{blue}{\textbf{[RM] #1}}}
\newcommand{\HW}[1]{\textcolor{green}{\textbf{[HW] #1}}} 		  
\else
\newcommand{\RM}[1]{}
\newcommand{\HW}[1]{} 	
\fi

\usepackage{csquotes}
\usepackage{verbatim}

\usepackage{epsfig}
\usepackage{color}

\usepackage{balance}  
\usepackage{algorithmic}
\usepackage{algorithm}
\usepackage{multirow}
\usepackage{graphicx}
\usepackage{amsmath}
\usepackage{amssymb}
\usepackage{amsfonts}
\usepackage{xspace}
\usepackage{url}
\usepackage{bbm}
\usepackage{stmaryrd}

\newcommand{\R}{\mathbb{R}}
\newcommand{\T}{\top}
\newcommand{\la}{\langle}
\newcommand{\ra}{\rangle}

\newcommand{\lam}{\lambda}

\newcommand{\Dt}{\Delta}

\newcommand{\conv}{{\rm conv}}

\newcommand{\al}{\alpha}

\newcommand{\ga}{\gamma}
\newcommand{\ep}{\epsilon}

\newcommand{\na}{\nabla}
\newcommand{\lt}{\left}
\newcommand{\rt}{\right}

\newcommand{\argmin}{\mathop{{\rm argmin}}}

\newcommand{\td}{\tilde}
\newcommand{\wtd}{\widetilde}

\newcommand{\Pb}{\mathbb{P}}

\newcommand{\supp}{\mathsf{supp}}

\newcommand{\cH}{\mathcal{H}}

\newcommand{\cJ}{\mathcal{J}}
\newcommand{\cK}{\mathcal{K}}

\newcommand{\cS}{\mathcal{S}}

\newcommand{\cL}{\mathcal{L}}

\newcommand{\bba}{{a}}
\newcommand{\bbb}{{b}}

\newcommand{\bbe}{{e}}

\newcommand{\bbp}{{p}}
\newcommand{\bbq}{{q}}

\newcommand{\bbu}{{u}}
\newcommand{\bbv}{{v}}
\newcommand{\bbw}{{w}}
\newcommand{\bbx}{{x}}
\newcommand{\bby}{{y}}
\newcommand{\bbz}{{z}}
\newcommand{\bbone}{{1}}
\newcommand{\bblam}{{\lambda}}
\newcommand{\bbmu}{{\mu}}

\usepackage{xcolor}

\newtheorem{theorem}{Theorem}[section]

\newtheorem{lemma}[theorem]{Lemma}

\def\endproof{\hfill $\Box$ \vskip 0.4cm}

\makeatletter
\@addtoreset{equation}{section}
\makeatother



\begin{document}

	\title{A Cyclic Coordinate Descent Method for Convex Optimization on Polytopes
		\footnote{This research is supported in part by a grant from the Office of Naval Research (N000142212665).}
	}
	\author{
		Rahul Mazumder\thanks{MIT Sloan School of Management, Operations Research Center and MIT Center for Statistics~({email: {rahulmaz@mit.edu}}).}
		\and
		Haoyue Wang\thanks{MIT Operations Research Center (email: haoyuew@mit.edu).}
	}

\date{}
\maketitle	

\begin{abstract}
		Coordinate descent algorithms are popular for huge-scale optimization problems due to their low cost per-iteration. Coordinate descent methods apply to problems where the constraint set is separable across coordinates. In this paper, we propose a new variant of the cyclic coordinate descent method that can handle polyhedral constraints provided that the polyhedral set does not have too many extreme points such as $\ell_1$-ball and the standard simplex. Loosely speaking, our proposed algorithm PolyCD, can be viewed as a {\emph{hybrid}} of cyclic coordinate descent and the Frank-Wolfe algorithms. 
 We prove that PolyCD has a $O(1/k)$ convergence rate for smooth convex objectives. Inspired by the away-step variant of Frank-Wolfe, we propose PolyCDwA, a variant of PolyCD with away steps which has a linear convergence rate when the loss function is smooth and strongly convex. Empirical studies demonstrate that PolyCDwA achieves strong computational performance for large-scale benchmark problems including $\ell_1$-constrained linear regression, $\ell_1$-constrained logistic regression and kernel density estimation. 
\end{abstract}

	\section{Introduction}
	\label{section: Introduction}

	The coordinate descent (CD) algorithm \cite{bertsekas1997nonlinear,wright2015coordinate} has a history that dates back to the origin of nonlinear optimization. A basic version of CD updates one coordinate at a time (to optimize the objective), while keeping other coordinates fixed.
	In the past ten years or so CD algorithms have been found to be very useful in various large-scale optimization problems arising in practice, particularly in statistics and machine learning \cite{friedman2010regularization,nesterov2012efficiency}. Due to its low per-iteration cost,
	in many large-scale applications (e.g. \cite{friedman2010regularization,mazumder2011sparsenet}), CD methods can outperform
	methods that require full gradient or Hessian evaluations \cite{wright2015coordinate,nesterov2012efficiency}. 
	For example, CD methods have state-of-the-art performance on large-scale $\ell_1$-regularized linear/logistic regression \cite{friedman2010regularization} and are widely implemented within packages such as glmnet and scikit-learn.

	When solving unconstrained convex optimization problems with smooth objectives, CD methods can be applied with convergence guarantees, see e.g.~\cite{luo1992convergence,luo1993error} for earlier results and~\cite{nesterov2012efficiency,beck2013convergence,hong2017iteration} for more recent results. However, when considering constrained convex optimization problems, 
	CD methods may not directly apply. For example, consider minimizing the two-dimensional function
	$f(x_1,x_2):= (x_1-2)^2 + (x_2-2)^2$ with the $\ell_{1}$-norm constraint $|x_1|+|x_2|\le 1$.
	If we initialize at $x_1=0,x_2=1$, then any CD step updating one coordinate at a time will violate the constraint. Therefore, a direct application of CD does not appear to converge for this example. 
	We encounter similar problems in the constrained form of the LASSO~\cite{tibshirani1996regression} while minimizing a least squares loss 
	with an $\ell_1$-ball constraint (for example). Interestingly, while the CD method is very popular and effective for an $\ell_1$-norm {\emph{penalized}} version of LASSO~\cite{friedman2010regularization}, this method does not apply to the original constrained LASSO version\footnote{There is a variant of CD method called 2-coordinate-descent method \cite{beck20142}, which can be applied to LASSO if one rewrites LASSO as an optimization on the standard simplex. See Section~\ref{section:literature-CD} for discussions on this approach.}.

	In this paper, we extend the scope of CD algorithms to consider a general problem of the form:
	\begin{equation}\label{problem-intro}
		\min ~~ f(\bbx) \quad \text{s.t.} ~ \bbx \in S
	\end{equation}
	where $f(\cdot)$ is a convex smooth (i.e., with Lipschitz-continuous gradient) function, and $S$ is a bounded polytope whose vertices are given by $\{\bbv^1,\ldots,\bbv^M\}$. 
	In this paper, we make the underlying assumption that the number of vertices $M$ is not too large. 
	Problem~\eqref{problem-intro} includes many important 
	instances---
	See Section~\ref{section: Applications} for applications in statistics and machine learning where such problems arise frequently. 
	For Problem~\eqref{problem-intro}, the classic CD method is not applicable, and 
	we propose a new variant of CD, namely a \textit{polyhedral coordinate descent (PolyCD)} method for this problem. 
	We present some intuition underlying PolyCD. 
	The vertices in the polytope $S$ correspond to the coordinates in an unconstrained problem. At every iteration of PolyCD, we pick a vertex and move in the direction toward vertices (from the current iterate). We cycle across the different vertices in a spirit similar to a cyclic CD (CCD) method.
	Section \ref{section: Polyhedral coordinate descent algorithm} presents a full description of PolyCD. 
	In addition to the similarity to CD methods, 
	PolyCD also has some similarities to the cyclic variant of the Frank-Wolfe (FW) method~\cite{frank1956algorithm,jaggi2013revisiting}. 
	Intuitively, the vanilla FW method has similarities with the greedy coordinate descent method, while our method has close parallels to cyclic coordinate descent. See Section~\ref{subsection:literature-FW} for additional discussions.

	We present an instantiation of Problem~\eqref{problem-intro} that poses computational challenges for large-scale problems where PolyCD seems to be quite promising. 
	Consider the following problem:
	\begin{equation}\label{problem2}
		\begin{aligned}
			\min_{\bbx} ~~  f(\bbx): = g(A \bbx) 
			\quad \text{s.t.} ~ x \in \Dt_d
		\end{aligned}
	\end{equation}
	where $g(\cdot)$ is a smooth function and $A\in \R^{n\times d}$ is a dense matrix, and $\Dt_d$ is the standard simplex $\Dt_d:= \{ \bbx\in \R^d~|~ \bbx\ge 0, ~ {1}_d^\T \bbx =1\}$. Problems of this form appear in, e.g., convex approximation \cite{zhang2003sequential}, core vector machines \cite{tsang2006generalized} and Adaboost \cite{freund1997decision}
	(see Section \ref{section: Applications} for more applications). Most existing methods for problem \eqref{problem2} require (at least) the evaluation of gradient in each iteration. Since $\na f(\bbx) = A^\T \na g(A \bbx)$, the cost for one gradient evaluation is typically $O(nd)$ -- this can be expensive when both $n$ and $d$ are large. In comparison, given the current iterate $\bbx^0$ and a vertex $\bbe_i$, 
	the one-dimensional optimization on the line segment with endpoints $\bbx^k$ and $\bbe_i$ only requires $O(n)$ operations, so the per-iteration cost of PolyCD is much smaller than most first-order methods. 
	Note that CCD methods have similar advantages when the objective function is separable, but due to the nature of the constraint, it cannot be directly applied to Problem~\eqref{problem2}.
	Empirically, we find that our algorithm typically requires only a few passes over the data, and can achieve a significant speedup over existing methods on many benchmark 
	problems---See Section~\ref{section: Experiments} for details.

	\subsection{Related literature}
	We provide an overview of work on FW and CD methods as they are related to our proposed approach.

	\subsubsection{Frank-Wolfe methods}\label{subsection:literature-FW}
	Frank-Wolfe method is a first-order method for smooth optimization on bounded constraint set proposed in the 1950s \cite{frank1956algorithm}. Recently it has received significant attention in the optimization and machine learning communities due to its projection-free nature. In each iteration, FW computes a linear subproblem over the constraint set.
	When the constraint set is a polytope, an optimal solution to the linear subproblem is at a vertex. 
	The FW algorithm 
	makes a move on the line segment joining the current point and this vertex. 
	Note that FW can also be applied to Problem~\eqref{problem-intro}, but requires evaluation of the full gradient in each iteration, which can be expensive for the applications we consider (e.g. $\ell_1$-constrained logistic regression). 
	
	There are some similarities between the FW and greedy CD 
	method~\cite{nutini2015coordinate,lu2018accelerating}---at every iteration, both methods compute the full gradient and choose a coordinate in a greedy fashion. 
	In the literature, there are roughly two other  types of CD methods: the cyclic CD (CCD) method and the randomized CD (RCD) method. Different CD methods have different operating characteristics under different settings -- see e.g. \cite{sun2019worst,nutini2015coordinate,gurbuzbalaban2017cyclic} for discussions. In particular, 
	since CCD and RCD do not require the evaluation of a full gradient in every step, their per-iteration cost is much smaller than the greedy CD. This difference can be significant when gradient evaluations are expensive. 
	Even though FW, as a counterpart of greedy CD, has been widely used for problem \eqref{problem-intro}, surprisingly, 
	there appears to be no counterpart to CCD in the context of 
	Problem~\eqref{problem-intro} (as far as we know). 
	In this paper, we seek to fill this gap in the literature.

	When the objective function is smooth and convex, FW has a convergence rate $O(1/t)$~\cite{jaggi2013revisiting}, where $t$ denotes the iteration index. 
	To improve the vanilla version of FW, Away-step Frank-Wolfe (AFW) method \cite{guelat1986some} and other variants \cite{lacoste2015global} have been proposed. These variants incur slightly more computational cost and memory in each iteration but are typically (overall) much faster than the vanilla FW method. When the objective function is strongly convex, and the constraint set is a polytope, AFW has a global linear convergence rate \cite{lacoste2015global,beck2017linearly}. Inspired by these variants of FW, we also propose a variant of PolyCD by introducing away steps (see Section \ref{section: Polyhedral coordinate descent with away steps} for details).

	\subsubsection{Coordinate descent methods}\label{section:literature-CD}
	While CD algorithms have been around for several years, the iteration complexity of CCD, especially the version with full minimization steps for smooth convex functions, 
	have been studied fairly recently~\cite{beck2013convergence,hong2017iteration,saha2013nonasymptotic}.
	Several existing analysis of CCD methods relies on an error-bound argument (see e.g. \cite{luo1993error,beck2013convergence,hong2017iteration}). 
	Although the worst-case iteration complexity of CCD given by existing literature is inferior to gradient descent (GD) methods ($\#$ cycles for CCD vs. $\#$ gradients steps for GD) by a factor dependent on the dimension $d$, the practical performance of CCD might be much better~\cite{friedman2010regularization,hazimeh2020fast}. 
	For example, CCD is used in glmnet~\cite{friedman2010regularization} for the unconstrained form of LASSO\footnote{That is, we penalize the $\ell_{1}$-norm of the regression coefficients instead of including it in the form of a bound constraint.} and appears to converge faster than the worst-case rate as suggested by theory.
	Similarly, in our experiments, our proposed algorithm appears to run faster than what is suggested by our theoretical guarantees.

	There is an interesting line of research \cite{platt1998sequential,keerthi2001improvements,beck20142,necoara2017random,reddi2014large,cristofari2019almost,cristofari2022active} that generalizes coordinate descent for the minimization of $f(x) + h(x)$ over constraints of the form: $Ax=b$, where $f(\cdot)$ is a smooth convex function, $h(\cdot)$ is a separable (possibly, non-smooth) convex function, and $A\in \R^{m\times d}$ with $m$ being very small. The minimization of $f(x)$ on the standard simplex is a specific example of this setting, with $A = 1_{d}^\T$, $b=1$, and $h(x) = \sum_{i=1}^d \iota(x_i)$, where $\iota(t) = 0$ for $t\ge 0$ and $\iota(t) = \infty$ for $t<0$. For this problem, 
	the algorithms in \cite{beck20142,necoara2017random,reddi2014large} take ``pairwise" coordinate descent steps where only a pair of coordinates are modified in each iteration---this is the so-called ``2-CD" method. The 
	2-CD method differs from our PolyCD approach, as we discuss below. First, the 2-CD method is very similar to a CD method and aims to modify only a small number of coordinates (depending on $m$) in each step. 
	Our method adopts a different intuition and moves towards a vertex in each step.
	Second, the theoretical guarantees for our approach and 2-CD methods are different. Both our method and 2-CD are known to have a sublinear rate $O(1/t)$ when the objective function is convex and smooth \cite{beck20142,reddi2014large}. When the function is also strongly convex, we prove a non-asymptotic linear rate of our algorithm (See Section~\ref{section: Polyhedral coordinate descent with away steps} for details). It seems that an asymptotic linear rate of convergence has been established for a recently proposed variant of the 2-CD method \cite{cristofari2019almost}\footnote{Note that \cite{cristofari2019almost} derive a non-asymptotic linear rate for an unbounded problem with a single linear constraint, which does not include the standard simplex (standard simplex is bounded).}.
	Finally, the numerical performance of our methods and pairwise CD can be quite different. 
	Section~\ref{subsection: exp: l1 constrained least squares} presents an empirical comparison of the 2-CD method in \cite{necoara2017random,reddi2014large} and our method for the LASSO.

	\subsection{Notations and preliminaries}
	
	Let $\|\cdot\|$ denote the Euclidean norm. For any two points $\bbx,\bby\in \R^d$, let $[\bbx,\bby]$ be the line segment with end points $\bbx$ and $\bby$. 
	Given integer $m>0$, let $\bbone_m$ be the vector in $\R^m$ with all coordinates being $1$. Let $e_1,\ldots,e_d$ be vectors in $\R^d$, where $e_i$ has its $i$-th coordinate being $1$ and all other coordinates being $0$. 
	Let $\mathbb{S}^{d\times d}_+$ be the set of positive semidefinite matrix in $\R^{d\times d }$. For an integer $M \geq 1$, we use the notation $[M]=\{1, \ldots , M\}$.
	Given a convex set $S\subseteq \R^d$ and $L>0$, a function $f(\cdot)$ is called \textit{$L$-smooth on $S$} if it is differentiable and 
	$\|\na f(\bbx) - \na f(\bby)\| \le L \|\bbx- \bby\|$
	for all $x,y\in S$.  
	Given $\mu>0$, 
	a function $f(\cdot)$ is called \textit{$\mu$-strongly-convex} on S if  $f(\bby) - f(\bbx) \ge \la \na f(\bbx), \bby- \bbx \ra + \frac{\mu}{2} \| \bby- \bbx\|^2$ 
	for all $ \bbx, \bby \in S$.

	The remainder of the paper is organized as follows. 
	In Section \ref{section: Polyhedral coordinate descent algorithm} we formally describe our PolyCD algorithm and prove its global $O(1/k)$ convergence rate under standard smoothness assumptions of $f(\cdot)$. In Section~\ref{section: Polyhedral coordinate descent with away steps}, we propose an improved version of PolyCD with away steps (denoted as PolyCDwA), and prove that it is globally linear convergent under the assumption that $f(\cdot)$ is strongly convex. In Section \ref{section: Applications}, we discuss applications of our framework, and in Section \ref{section: Experiments}, we present the numerical performance of our proposed methods and comparisons with existing methods.

	\section{Polyhedral coordinate descent algorithm (PolyCD)}\label{section: Polyhedral coordinate descent algorithm}
	In this section, we formally present our proposed polyhedral coordinate descent algorithm (Algorithm~\ref{alg: PolyCD}) for Problem~\eqref{problem-intro} and prove its convergence rate for smooth and convex loss functions. Recall that $\{\bbv^1, \ldots, \bbv^M\}$ is the set of extreme points of $S$.

	\begin{algorithm}[H]
		\caption{Polyhedral coordinate descent method (PolyCD) for problem \eqref{problem-intro}}
		\label{alg: PolyCD}
		\begin{algorithmic}
			\STATE Start from $\bbx^{0,0} \in S$. 
			
			\STATE For $t=0,1,2,....$, 
			\STATE \quad For $i=1,2,....,M$, update:
			\vspace{-0.2cm}
			\begin{equation}\label{PolyCD: update}
				\vspace{-0.2cm}
				\bbx^{t, i} =  \bbx^{t,i-1} + \al_{t,i}  ( \bbv^i - \bbx^{t,i-1})  
			\end{equation}
			\qquad with $\al_{t,i} $ chosen in $[0,1]$. 
			\STATE \quad 
			Set $\bbx^{t+1, 0} = \bbx^{t, M}$. 
		\end{algorithmic}
	\end{algorithm}
	In the outer iteration $t$, Algorithm \ref{alg: PolyCD} sequentially performs the update~\eqref{PolyCD: update} across 
	$i \in [M]$ where $i$ indexes the extreme points of $S$. In particular, for every $i$, update~\eqref{PolyCD: update} moves the current solution towards the $i$-th extreme point $v_{i}$ with corresponding step-size $\al_{t,i}$.
	The step size $\al_{t,i}$ 
	can be chosen by different rules. For example, it can be computed by an exact line search over the line segment $[\bbx^{t,i-1}, \bbv^i]$: 
	\begin{equation}\label{stepsize: line-search-rule}
		\al_{t,i} \in \argmin_{\al\in [0,1]}  f( \bbx^{t,i-1} + \al ( \bbv^{i} - \bbx^{t,i-1}) )  .
	\end{equation}
	Update \eqref{stepsize: line-search-rule} 
	can be computed in closed form for some special cases of $f(\cdot)$ (e.g. quadratic functions). 
	For a general function $f(\cdot)$, the line-search update \eqref{stepsize: line-search-rule} may require multiple function evaluations on the line segment $[\bbx^{t,i-1}, \bbv^i]$, which can be expensive for large-scale problems. 
	In such cases, for an $L$-smooth function $f(\cdot)$, we consider an alternative step size rule by performing a proximal gradient step on the line segment 
	$[\bbx^{t,i-1}, \bbv^i]$: 
	\begin{equation}\label{stepsize: 1d-gradient-rule}
		\al_{t,i} \in \argmin_{\al\in [0,1]} \Big\{   \al \la  \na f(\bbx^{t,i-1}),  \bbv^{i} - \bbx^{t,i-1}\ra + \frac{L\al^2}{2} \|  \bbv^{i} - \bbx^{t,i-1} \|^2   \Big\}.
	\end{equation}
	To compute the 
	update \eqref{stepsize: 1d-gradient-rule}, instead of computing the full gradient $\na f(\bbx^{t,i-1})$, one only needs to compute the value $\la  \na f(\bbx^{t,i-1}),  \bbv^{i} - \bbx^{t,i-1}\ra$. In many applications (see Section \ref{section: Applications}), computing the latter can be computationally friendlier than computing the full gradient. Hence the per-iteration cost of PolyCD can be much lower than first-order methods which require evaluations of full gradients.

	\subsection{Convergence guarantees}
	We state and prove the convergence rate of PolyCD (Algorithm \ref{alg: PolyCD}) for both step-size rules \eqref{stepsize: line-search-rule} and \eqref{stepsize: 1d-gradient-rule}. 
	In the following, for the iterations $\{\bbx^{t,i}\}$ generated by Algorithm \ref{alg: PolyCD}, we use the notation $\bbx^t := \bbx^{t,0}$ for all $t\ge 0$. Denote $D:= \sup_{\bbx,\bby\in S} \| \bbx-\bby\|$. Let $\bbx^*$ be an optimal solution of \eqref{problem-intro}, and $f^* = f(\bbx^*)$.

	\begin{theorem}\label{theorem: PolyCD}
		(Sublinear rate)
		Suppose $f(\cdot)$ is convex and $L$-smooth on $S$. 
		Let $\{\bbx^{t,i}\}_{t \ge 0, 0\le  i  \le M}$ be the sequence of iterates generated by 
		Algorithm~\ref{alg: PolyCD}. Then the following holds true:
		
		(1) If exact line search steps \eqref{stepsize: line-search-rule} are used, then for all $t\ge 1$, 
		\begin{equation}\label{PolyCD converge 1}
			f(\bbx^t) - f^* ~ \le ~ \frac{\max\{ f(\bbx^1) - f^*, 4MLD^2 \}}{t} .
		\end{equation}
		
		(2) If one-dimensional gradient steps \eqref{stepsize: 1d-gradient-rule} are used, then for all $t\ge 1$, 
		\begin{equation}\label{PolyCD converge 2}
			f(\bbx^t) - f^* ~ \le ~ \frac{\max\{ f(\bbx^1) - f^*, 16MLD^2 \}}{t}.
		\end{equation}
	\end{theorem}
	
	Theorem \ref{theorem: PolyCD} states that for iterates $\bbx^t(= \bbx^{t,0})$ generated by Algorithm \ref{alg: PolyCD}, the optimality gap
	$f(\bbx^t) - f^*$ converges to $0$ with the rate $O(1/t)$
	for both step-size rules \eqref{stepsize: line-search-rule} and \eqref{stepsize: 1d-gradient-rule}. 
	The upper bounds on $f( \bbx^t) - f^*$ in \eqref{PolyCD converge 1} and \eqref{PolyCD converge 2} depend on the smoothness parameter $L$, the diameter $D$ of the constraint $S$, and the number of extreme points $M$. 
	As $M$ can be large in some applications, so can the upper bounds. For example, for least squares with an $\ell_1$-norm constraint, the number of features $p$ (and hence, $M = 2p$) can be large. Indeed, the CCD method has a similar convergence rate\footnote{We omit other constants in the rate and only show the dependence on $p$ and $t$} $O(p/t)$ where the outer iterations are indexed by $t$~\cite{beck2013convergence,hong2017iteration}.
	Unlike the rates in~\eqref{PolyCD converge 1} and \eqref{PolyCD converge 2}, the rate of the FW method does not have the $M$ dependence. This is not surprising as the FW method requires computing the full gradient at every step, which differs from our setting. 
	Interestingly, a similar gap in worst-case computational guarantees exists between greedy CD and cyclic CD methods.

	\section{PolyCD with away steps (PolyCDwA)}\label{section: Polyhedral coordinate descent with away steps}
	The basic version of PolyCD (Algorithm \ref{alg: PolyCD}) may experience slow convergence as the iterations progress. To gather intuition, consider the problem of least squares regression
	with an $\ell_1$ norm constraint on the regression coefficients. 
	As PolyCD updates the coordinates in a cyclic order, it is possible that in the first few iterations, PolyCD takes large steps toward vertices that should not be in the final support of the optimal solution. Since PolyCD only moves toward vertices, it may take a long time to ``offset" the first few ``bad" steps. This leads to a slow convergence of PolyCD, which is also observed in our numerical experiments (Section~\ref{subsubsection: linear regression Comparison between PolyCD and PolyCDwA}). 
	
	To fix this problem, we propose an improved version of PolyCD by extending the line segment $[x^{t,i-1},v^i]$ to a larger line segment (while remaining within the constraint $S$) and making the update on this extended line segment. 
	In other words, this extended line segment allows for backward steps that move away from a given vertex. 
	We call our proposed algorithm \textit{Polyhedral Coordinate Descent method with Away steps} (abbreviated as PolyCDwA), summarized in 
	Algorithm~\ref{alg: PolyCDwA} below. Note that the away steps in Algorithm~\ref{alg: PolyCDwA} resemble the away steps in the away-step variant of FW (AFW). But unlike AFW,
	which computes the full gradient to choose an away step direction, PolyCDwA
	maintains the cyclic nature of PolyCD and the low per-iteration cost for each coordinate update.

	\begin{algorithm}[H]
		\caption{Polyhedral Coordinate Descent method with Away steps (PolyCDwA)}
		\label{alg: PolyCDwA}
		\begin{algorithmic}
			\STATE Start from $\bbx^{0,0} \in S$, with $\bbx^{0,0} = \sum_{i=1}^M \lam^{0,0}_i \bbv^{i}$ satisfying $\lam_i^{0,0} \ge 0$ 
			and $\sum_{i=1}^M \lam_i^{0,0} = 1$. 
			
			\STATE For $t=0,1,2,....$, 
			
			\STATE \quad For $i = 1, 2, ...,M$, 
			let \vspace{-0.2cm}
			\begin{equation}\label{PolyCDwA def gamma}
				\vspace{-0.2cm}
				\ga_{t,i} = \lam_{i}^{t,i-1} / (1- \lam_{i}^{t,i-1}),  
			\end{equation}
			\qquad and $\ga_{t,i} = \infty$ if $\lam_{i}^{t,i-1} = 1$. 
			Update \vspace{-0.2cm}
			\begin{equation}\label{PolyCDwA update x}
				\vspace{-0.2cm}
				\bbx^{t, i} =  \bbx^{t,i-1} + \al_{t,i}  (\bbv^{i} - \bbx^{t,i-1}) 
			\end{equation}
			\qquad with $\al_{t,i} $ chosen in $[-\ga_{t,i},1]$, and
			\vspace{-0.2cm}
			\begin{equation}\label{PolyCDwA update lam}
				\vspace{-0.2cm}
				\begin{aligned}
					& \lam_{i}^{t,i} = (1-\al_{t,i}) \lam_i^{t,i-1} + \al_{t,i} \\
					& 	\lam_{j}^{t,i} = (1-\al_{t,i}) \lam_j^{t,i-1} ~~ \forall ~ j\in [M] \setminus \{i\}. 
				\end{aligned}
			\end{equation}
			
			\vspace{0.2cm}
			\STATE \quad Let $\bbx^{t+1, 0} = \bbx^{t, M}$. 
		\end{algorithmic}
	\end{algorithm}
	In Algorithm~\ref{alg: PolyCDwA}, the indices $t$ and $i$ stand for outer and inner iterations, respectively (this is similar to our notation for Algorithm~\ref{alg: PolyCD}).
	Algorithm \ref{alg: PolyCDwA} maintains a decomposition of the current iterate in terms of extreme points:
	$\bbx^{t,i} = \sum_{j=1}^{M} \lam_{j}^{t,i} \bbv^i$, where ${\lam}^{t,i} := [\lam^{t,i}_1, ..., \lam^{t,i}_M]^\T \in \Dt_M$ (this can be verified by the fact that $\bblam^{0,0} \in \Dt_M$ and the updates in \eqref{PolyCDwA update lam}). 
	Making use of the decomposition
	$\bbx^{t,i-1} = \sum_{j=1}^{M} \lam_{j}^{t,i-1} \bbv^i$, 
	we first compute a value $\ga_{t,i}$ defined in \eqref{PolyCDwA def gamma}. The value $\ga_{t,i}$ is the largest step size that one can take to move away from the extreme point $\bbv^i$  while still remaining in the constraint set $S$. 
	With $\ga_{t,i}$ at hand, we perform the update in \eqref{PolyCDwA update x}. This update is similar to the update~\eqref{PolyCD: update} in vanilla PolyCD, but we allow $\al_{t,i}$ to take negative values that are larger than $-\ga_{t,i}$. 
	In particular, when $\al_{t,i}<0$, the iterate moves away from $\bbv^i$ along the line joining $\bbx^{t,i-1}$ and $\bbv^i$---this is what we refer to as an ``away step". 
	After making the update in \eqref{PolyCDwA update x}, we accordingly compute $\lam^{t,i}$ such that the representation $\bbx^{t,i} = \sum_{j=1}^{M} \lam_{j}^{t,i} \bbv^i$ holds. 
	
	Similar to Algorithm \ref{alg: PolyCD}, we consider two ways to select the step size $\al_{t,i}$. One approach is to use exact line-search
	\begin{equation}\label{stepsize: line-search-rule-away}
		\al_{t,i} \in \argmin_{\al\in [-\ga_{t,i},1]} \lt\{  f( \bbx^{t,i-1} + \al ( \bbv^{t,i} - \bbx^{t,i-1}) ) \rt\} 
	\end{equation}
	and another is to consider one-dimensional (proximal) gradient steps
	\begin{equation}\label{stepsize: 1D-gradient-rule-away}
		\al_{t,i} \in \argmin_{\al\in [-\ga_{t,i},1]} \Big\{  \al \la  \na f(\bbx^{t,i-1}), \bbv^{t,i} - \bbx^{t,i-1}\ra 
		+ \frac{L\al^2}{2} \|\bbv^{t,i} - \bbx^{t,i-1} \|^2   \Big\}.
	\end{equation}
	Note that the updates in \eqref{stepsize: line-search-rule-away} and \eqref{stepsize: 1D-gradient-rule-away} differ from the updates in \eqref{stepsize: line-search-rule} and \eqref{stepsize: 1d-gradient-rule} in the range of $\al$ values considered. 
	Compared to PolyCD (i.e., Algorithm \ref{alg: PolyCD}), PolyCDwA (i.e., Algorithm \ref{alg: PolyCDwA}) 
	incurs the additional overhead of maintaining the decomposition in terms of the extreme points. In particular, 
	we need to maintain a weight vector $\bblam^{t,i}$, which requires $O(M)$ memory. In addition, there is  a $O(M)$ cost in updating $\bblam^{t,i}$ by \eqref{PolyCDwA update lam}.
	For most of the applications we considered (see Section \ref{section: Applications}), the per-iteration costs of PolyCD and PolyCDwA are comparable. Still, the overall empirical performance of PolyCDwA appears to be much better (see Section \ref{section: Experiments}).

	\subsection{Computational guarantees}
	Below we present computational guarantees for PolyCDwA. 
	Recall that we use the notation: $\bbx^t= \bbx^{t,0}$ for all $t\ge 0$; $\bbx^*$ denotes an optimal solution to~\eqref{problem-intro}, and $f^* = f(\bbx^*)$. 
	
	Theorem~\ref{theorem: PolyCDwA weak convexity} shows that when $f$ is convex and smooth, PolyCDwA has the same convergence rate as PolyCD. The proof of Theorem \ref{theorem: PolyCDwA weak convexity} is (almost) the same as the proof of Theorem \ref{theorem: PolyCD} and hence omitted for simplicity. 
	\begin{theorem}\label{theorem: PolyCDwA weak convexity}
		(Convex, sublinear rate)
		Suppose $f(\cdot)$ is convex and $L$-smooth on $S$. 
		Let $\{\bbx^{t,i}\}_{t \ge 0, 0\le  i  \le M}$ be the sequence generated by 
		Algorithm~\ref{alg: PolyCDwA}. 
		
		(1) If line search steps \eqref{stepsize: line-search-rule-away} are used, then for all $t\ge 1$, 
		\begin{equation}\label{PolyCDwA weak convexity converge 1}
			f(\bbx^t) - f^* ~ \le ~ \frac{\max\{ f(\bbx^1) - f^*, 4MLD^2 \}}{t} .
		\end{equation}
		
		(2) If one-dimensional gradient steps \eqref{stepsize:  1D-gradient-rule-away} are used, then for all $t\ge 1$,
		\begin{equation}\label{PolyCDwA weak convexity converge 2}
			f(\bbx^t) - f^* ~ \le ~ \frac{\max\{ f(\bbx^1) - f^*, 16MLD^2 \}}{t} .
		\end{equation}
	\end{theorem}
	
	Below we explore the convergence of PolyCDwA when $f(\cdot)$ is strongly convex. 
	First, we define the facial distance $\psi_S$ \cite{pena2019polytope} of the polyhedral set $S$ as follows:
	\begin{equation}\label{facial-distance-defn}
		\psi_S := \min_{\mbox{$\scriptsize{
					\begin{array}{c}
						F \in \text{faces}(S), \\
						\emptyset \neq F \neq S.
					\end{array}
				} $}}
		\text{dist}(F , \text{conv}(V(S) \backslash F)) \ ,
	\end{equation}
	where $V(S)$ denotes the set of all vertices of $S$; $\text{faces}(S)$ denotes the set of all faces of $S$.
	Then we have the following theorem on the convergence of PolyCDwA. 
	
	\begin{theorem}\label{theorem: PolyCDwA strong convex}
		(Strongly convex, linear rate)
		Suppose $f(\cdot)$ is $L$-smooth and $\mu$-strongly convex on $S$. 
		Let $\{\bbx^{t,i}\}_{t \ge 0, 0\le  i  \le M}$ be the sequence generated by Algorithm~\ref{alg: PolyCDwA}. 
		
		(1) If line search steps \eqref{stepsize: line-search-rule-away} are used, defining $G:=  1+ 9MLD^2 /(\mu \psi_S^2)   $, then 
		\begin{equation}
			f( \bbx^{t}) - f^* \le \Big(\frac{G}{1+G} \Big)^t (f( \bbx^0) - f^*) \quad \forall ~ t\ge 0.
		\end{equation}
		
		(2) If 1D gradient steps \eqref{stepsize: 1D-gradient-rule-away} are used, 
		defining $G':=  2 + 16MLD^2/(\mu \psi_S^2)$, we have: 
		\begin{equation}
			f( \bbx^{t}) - f^* \le \Big(\frac{G'}{1+G'} \Big)^t (f( \bbx^0) - f^*) \quad \forall ~ t\ge 0.
		\end{equation}
	\end{theorem}
	
	Theorem \ref{theorem: PolyCDwA strong convex} shows the linear convergence of PolyCDwA under the strong convexity assumption on $f(\cdot)$ for both step-size rules \eqref{stepsize: line-search-rule-away} and \eqref{stepsize: 1D-gradient-rule-away}. In addition to the dependence on $M$, $L$, $D$ and $\mu$, the convergence rate parameters $G$ and $G'$ also depend on the geometric constant $\psi_S$. 
	This is similar to the convergence rates for AFW \cite{lacoste2015global} since both analyses of PolyCDwA and AFW make use of a condition number of the objective function relative to the constraint set \cite{pena2019polytope}. Nevertheless, the overall proof techniques for AFW and for PolyCDwA are significantly different (see Section \ref{section: proof of theorem PolyCDwA strong convex} for details).

	Note that in both Theorems \ref{theorem: PolyCDwA weak convexity} and \ref{theorem: PolyCDwA strong convex}, the parameters in the upper bounds of $f(\bbx^t) - f^*$ depend on $M$. This appears to suggest that when $M$ is large the convergence of PolyCDwA is slow. However, in our numerical experiments, we empirically observed that the performance of PolyCDwA can be much better than the worst-case convergence rates given by Theorems \ref{theorem: PolyCDwA weak convexity} and \ref{theorem: PolyCDwA strong convex} (see Section \ref{section: Experiments}). 
	Finally, we note that the upper bounds in Theorems \ref{theorem: PolyCDwA weak convexity} and \ref{theorem: PolyCDwA strong convex} can be improved under some special assumptions. 
	For example, if all the iterations have small support $\cJ \subseteq [M]$, i.e., $\lam_{j}^{t,i} = 0$ for all $t\ge 0$ for $i\in [M]$ and $j\in [M] \setminus \cJ$,  
	then by a simple modification of the proof, the parameter $M$ in the upper bounds in Theorems \ref{theorem: PolyCDwA weak convexity} and \ref{theorem: PolyCDwA strong convex} can be replaced by a much smaller number $|\cJ|$.

	\section{Applications}\label{section: Applications}
	In this section, we present a few instantiations of Problem~\eqref{problem-intro} where PolyCD and PolyCDwA can be applied.

	\subsection{Optimization on the unit simplex}\label{subsection:application-opt-on-unit-simplex}
	The standard $d$-dimensional simplex $\Dt_d =  \{\bbx\in \R^d~|~ \bbx\ge 0, ~ \bbone_d^\T \bbx =1\}$ is a polytope whose set of vertices are given by $\{\bbe_1, \ldots , \bbe_d\}$. 
	Many applications arising in statistics, machine learning, computational geometry, and related fields can be formulated as a convex optimization on the standard simplex, 
	including core vector machines \cite{tsang2006generalized}, Adaboost \cite{freund1997decision,zhang2003sequential}, mixture density estimation \cite{li1999mixture}, minimum enclosing ball \cite{yildirim2008two}, and $L_p$ regression \cite{zhang2003sequential}.
	See \cite{clarkson2010coresets} for a survey on related applications. 
	
	To illustrate the computational cost of PolyCD and PolyCDwA, consider a loss function with a finite-sum structure:
	\begin{equation}\label{problem: ERM simplex}
		\min_{\bbx\in \R^d} ~~  f(\bbx):=  \sum_{j=1}^n h(\bba_j^\T \bbx) ~~~ \text{s.t.} ~~~ \bbx\in \Dt_d
	\end{equation}
	where $h(\cdot)$ is a smooth convex function on $\R$, and $\bba_i \in \R^d$ for all $i\in [n]$. Many concrete problems arising in practice~\cite{clarkson2010coresets} can be written in the form of \eqref{problem: ERM simplex}. In each iteration of PolyCD, we need to compute the
	step size $\al_{t,i}$. When $h(s) = s^2$, the line-search steps \eqref{stepsize: line-search-rule} can be calculated in closed form with $O(n+d)$ operations. More generally, when $h(\cdot)$ is a $L'$-smooth function with some $L'>0$, one can use the one-dimensional gradient steps \eqref{stepsize: 1d-gradient-rule}. To this end,
	it suffices to calculate two values
	\begin{equation}
		b_{t,i}:=  \la  \na f(\bbx^{t,i-1}),  \bbe_i - \bbx^{t,i-1}\ra , ~~
		c_{t,i}:= \|  \bbe_{i} - \bbx^{t,i-1} \|^2 \ .
		\nonumber
	\end{equation}
	It is easy to check that $c_{t,i}$ can be computed within $O(d)$ operations. 
	Denote $A:= [\bba_1, ...., \bba_n]^\T \in \R^{n\times d}$, then we have 
	$\na f(\bbx) = A^\T \xi$ with $\xi = [h'(\bba_1^\T \bbx), ...., h'(\bba_n^\T \bbx)]^\T \in \R^n$. If we keep a copy of $A\bbx^{t,i-1}$ in memory and update it in each iteration, then both $\la \na f(\bbx^{t,i-1}), \bbe_i \ra$ and $ \la \na f(\bbx^{t,i-1}), \bbx^{t,i-1} \ra$ (and hence $b_{t,i}$) can be computed with a cost of $O(n+d)$ operations. Therefore, when using the one-dimensional gradient steps~\eqref{stepsize: 1d-gradient-rule}, the per-iteration cost of PolyCD is $O(n+d)$. For the implementation of PolyCDwA, the only additional cost (over PolyCD) is maintaining and updating the weights $\bblam^{t,i}$, which takes $O(d)$ operations. So the per-iteration cost of PolyCDwA is also $O(n+d)$. 
	As a side note, for methods that require full gradient evaluations (e.g. the proximal gradient method), the per-iteration cost is $O(nd)$, which is much larger than both PolyCD and PolyCDwA. 
	
	\subsection{Optimization on the $\ell_1$-norm ball}
	The $\ell_1$-norm constraint or penalty commonly arises 
	in machine learning when sparse solutions are sought via convex optimization~
	\cite{hastie2015statistical}. 
	For example, the $\ell_1$-constrained least squares (aka LASSO) solves the problem
	\begin{equation}\label{problem: least squares l1 constraints}
		\min_{\bbx\in \R^d} ~~ f(\bbx):= \| A\bbx - \bbb\|^2 ~~\text{s.t.} ~~~ \|\bbx\|_1 \le C
	\end{equation}
	where $A \in \R^{n\times d}$, 
	$\bbb \in \R^n$, and $C>0$.  
	The constraint set $B_C:= \{\bbx\in \R^d~|~ \|\bbx\|_1 \le C\}$ has vertices $\{r\bbe_1, -r\bbe_1, \ldots, r\bbe_d, -r\bbe_d\}$. 
	The per-iteration costs of PolyCD and PolyCDwA for Problem~\ref{problem: least squares l1 constraints} are similar to that discussed in Section~\ref{subsection:application-opt-on-unit-simplex}.

	In the literature, the $\ell_1$-penalized version of LASSO seems to be more popular than the $\ell_1$-constrained version. The former appears to be computationally more appealing than the latter.
	In particular, for $\ell_1$-penalized least squares, cyclic CD algorithms \cite{friedman2010regularization} are known to be quite efficient and are suitable for large-scale problems. In contrast, 
	for Problem~\eqref{problem: least squares l1 constraints},
	there are no CD algorithms with similar efficiency. 
	In Sections~\ref{subsection: exp: l1 constrained least squares} and \ref{subsection: exp: l1 constrained logistic regression}, we present numerical experiments showing that the $\ell_1$-constrained version can be solved very efficiently with PolyCDwA.  In particular, we observe that the number of outer loops (i.e., the number $t$) needed by PolyCDwA is small and stable ($10\sim 50$) even for large-scale problems.

	\subsection{Other applications involving structured polytopes}
	Both PolyCD and PolyCDwA can be used in other problems arising in statistics.
	For example, for the estimation of sparse graphical 
	models~\cite{hastie2015statistical}, the constraint set is given by an $\ell_1$-ball (under symmetry constraints)---the vertices of this polytope can be enumerated. 
	Another family of examples arises in shape-restricted density estimation \cite{turnbull2014unimodal,wang2021nonparametric}. For one-dimensional density estimation, using Bernstein polynomial bases, several shape restrictions (e.g. monotonicity, concavity, unimodality) can be translated to corresponding shape restrictions on the weights of the bases elements---these shape restrictions are given by a polytope---See \cite{wang2021nonparametric} for details.

	\subsection{Subproblems in fully-corrective FW method}
	The fully-corrective Frank-Wolfe method (FCFW) \cite{lacoste2015global} is a variant of the FW method that is often used to improve the performance of the vanilla FW method. 
	In each iteration (say iteration $t$), FCFW maintains a finite set $\cS^{(t)}$ of extreme points such that the current iteration can be expressed as a linear combination of points in $\cS^{(t)}$. Then FCFW solves a subproblem that minimizes the objective value on the set $\conv(\cS^{(t)})$. Since the vertices of $\conv(\cS^{(t)})$ are already given by $\cS^{(t)}$ and the number of points in $\cS^{(t)}$ is usually not too large, PolyCD and PolyCDwA can be used to solve this subproblem. 
	A similar argument also holds for some other variants of FW methods, e.g. \cite{von1977simplicial,freund2017extended}.

	\section{Experiments}\label{section: Experiments}
	We present numerical experiments of PolyCD, PolyCDwA. We focus on three problems: (1) $\ell_1$-constrained least squares; (2) $\ell_1$-constrained logistic regression; (3) kernel density estimation, and show the computational results for these three applications in Sections \ref{subsection: exp: l1 constrained least squares}, \ref{subsection: exp: l1 constrained logistic regression} and \ref{subsection: exp: kernal density estimation} respectively. 
	Our code is written in Julia 1.2.0. All the computations were performed on the MIT engaging cluster with the assignment of $1$ CPU and 16GB RAM. 
	
	\subsection{$\ell_1$ constrained linear regression}\label{subsection: exp: l1 constrained least squares}
	Consider Problem \eqref{problem: least squares l1 constraints}
	where $A = [\bba_1, \bba_2, \dots, \bba_n]^\T \in \R^{n\times d}$ with each $\bba_i \in \R^d$; 
	$\bbb =[b_1,b_2,\dots, b_n]^\T \in \R^n$; and $C>0$.  
	The data is generated from the underlying model:
	\begin{equation}\label{data generation: least squares}
		b_i = \bba_i^\T \bbx^* + \ep_i, \quad i\in [n].
	\end{equation}
	Above, $\bba_1$, \ldots, $\bba_n$ are iid draws from a multivariate Gaussian distribution $N({0}_d, \Sigma)$, where the covariance matrix $\Sigma \in \R^{d\times d}$ has diagonal entries equal to $1$ and off-diagonal entries $0.1$. 
	The underlying coefficient $\bbx^*$ is a sparse binary vector with $\| \bbx^*\|_0=r$.
	The errors $\ep_i$ ($i\in [n]$) are iid from $N({0},\sigma^2)$ with some $\sigma>0$, and are independent of $A$. 
	We denote the \textit{Signal-to-Noise Ratio (SNR)} as the value $SNR:= \|A\bbx^*\|^2/(n\sigma^2)$.

	\subsubsection{Comparing PolyCD and PolyCDwA}\label{subsubsection: linear regression Comparison between PolyCD and PolyCDwA}
	We present an example comparing the performances of PolyCD and PolyCDwA on 
	Problem~\eqref{problem: least squares l1 constraints}. Data $A$ and $\bbb$ are generated with $n = d = 1000$, $r = 50$, $SNR = 10$, and we set $ C = \|\bbx^*\|_1$ to be the $\ell_1$-norm of the underlying true signal. 
	We also run MOSEK~\cite{andersen2000mosek} to get an estimate of the optimal value, which we denote by $f^*$.
	We measure the progress of PolyCD and PolyCDwA by the \textit{(relative) optimality gap}, which is defined as $(f(\bbx^t) - f^*)/\max\{|f^*|,1\}$ where $t$ denotes the outer iteration counter. 
	PolyCD and PolyCDwA are used with the line-search steps \eqref{stepsize: 1d-gradient-rule} and \eqref{stepsize: line-search-rule-away} respectively.

	\begin{figure}[h]
		\centering
		\includegraphics[width=0.5\textwidth]{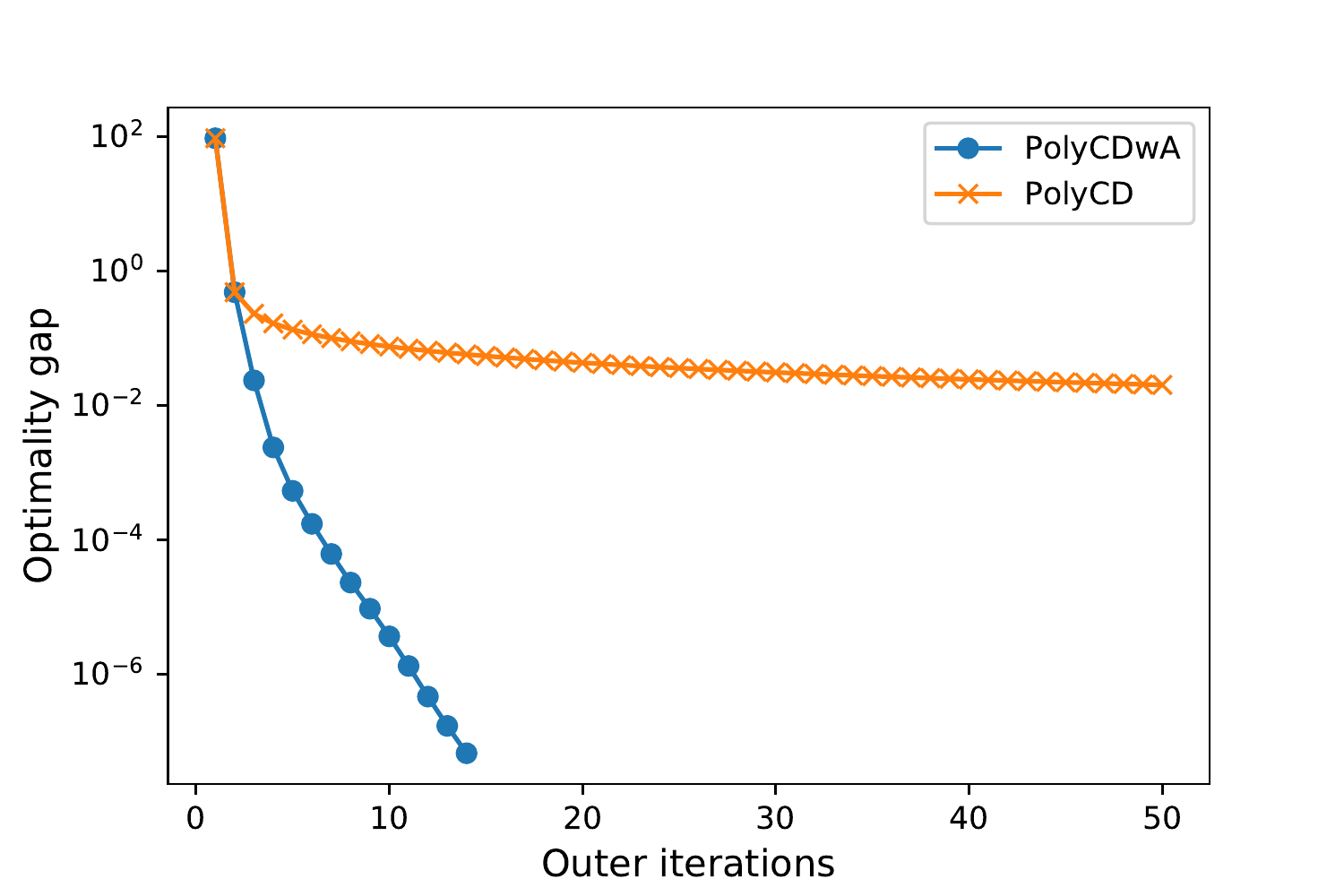}
		\caption{Comparison of PolyCD and PolyCDwA on $\ell_1$-constrained least squares.}
		\label{fig: compare PolyCD and PolyCDwA}
	\end{figure}
	
	Figure \ref{fig: compare PolyCD and PolyCDwA} presents the computational performances of PolyCD and PolyCDwA in this setting. It can be seen that the convergence of PolyCD is slow -- it cannot find a solution with an optimality gap less than $10^{-2}$ within $50$ outer iterations. This sublinear convergence performance is consistent with the theoretical results in Theorem \ref{theorem: PolyCD}. In contrast, PolyCDwA converges much faster and finds a solution with an optimality gap around $10^{-6} \sim 10^{-7}$ within 15 outer iterations. The convergence rate of PolyCDwA seems to be linear (as a reference, see Theorem \ref{theorem: PolyCDwA strong convex}). 
	
	\subsubsection{Comparison of PolyCDwA with benchmarks}\label{subsubsection: linear regression Comparison with benchmarks}
	We compare PolyCDwA with several state-of-the-art methods including: 
	2-CD~\cite{reddi2014large,necoara2017random}, AFW~\cite{lacoste2015global}, {FISTA}~\cite{beck2009fast}, a popular package~{StrOpt}\footnote{StrOpt stands for the Julia package StructuredOptimization.jl. It is a toolbox for structured optimization problems based on first-order methods. See  \url{https://github.com/JuliaFirstOrder/StructuredOptimization.jl} for details.}
	and the commercial solver MOSEK~\cite{andersen2000mosek}. 
	See Appendix~\ref{app-subsubsection: linear regression Comparison with benchmarks} for implementation details of these methods.

	We run PolyCDwA with a maximum of $100$ outer iterations and terminate it earlier if the relative improvement across two successive outer loops is less than $10^{-8}$. We let $f^*$ be the best objective value found across all algorithms.
	
	If $\hat f$ is the objective value computed by an algorithm, then the 
	\textit{optimality gap} of this algorithm is defined as $(\hat f - f^*) / \max\{|f^*|, 1\}$. 
	Table~\ref{table: linear regression compare with benchmarks} reports the runtimes (in seconds) and the optimality gaps of the algorithms
	on examples with $SNR=1$ and different values of $n$, $d$ and $r$. For StrOpt and MOSEK, only runtimes are reported. 
	We use different values of the regularization parameter $C$, each with different levels of sparsity---the corresponding number of nonzeros in the solution as computed by PolyCDwA are reported in the last column of Table~\ref{table: linear regression compare with benchmarks}. 
	The symbol ``x" stands for the instances where MOSEK runs out of memory. 
	All the reported numbers are the average of $5$ independent experiments. 
	As shown in Table~\ref{table: linear regression compare with benchmarks}, 
	PolyCDwA can find a high-accuracy solution (with an optimality gap less than $3\times 10^{-9}$) much more efficiently than other methods. In particular, PolyCDwA improves (in runtime) over the best of other methods by a factor of around $20\sim 100$.

	\begin{table}[h]
		\centering
		\caption{Comparison of PolyCDwA and benchmarks on $\ell_1$-constrained linear regression. Here, ``gap" stands for the optimality gap of an algorithm (see text for definition); ``Nonzeros" stands for the number of nonzero coordinates in the solution (as computed by PolyCDwA).}
		\scalebox{0.65}{
			\begin{tabular}{|c|cc|cccccccc|c|}
				\hline
				\multirow{2}{*}{}                                                             & \multicolumn{2}{c|}{PolyCDwA} & \multicolumn{2}{c}{2-CD} & \multicolumn{2}{c}{AFW} & \multicolumn{2}{c}{FISTA} & StrOpt & MOSEK & \multirow{2}{*}{Nonzeros} \\
				& time         & gap            & time       & gap         & time       & gap        & time       & gap          & time   & time  &                           \\ \hline
				\multirow{3}{*}{\begin{tabular}[c]{@{}c@{}}n=5K, d=5K,\\ r=500\end{tabular}}  & 0.6          & 2.0e-11        & 81.5       & 2.3e-03     & 15.6       & 2.8e-10    & 57.6       & 2.3e-08      & 20.9   & 172.9 & 233.4                     \\
				& 0.7          & 3.4e-12        & 99.9       & 1.4e-03     & 26.8       & 9.3e-10    & 67.2       & 4.9e-08      & 22.9   & 207.4 & 365.6                     \\
				& 0.9          & 3.7e-10        & 82.5       & 1.3e-03     & 34.7       & 1.6e-09    & 67.0       & 7.2e-08      & 24.3   & 197.8 & 451.2                     \\ \hline
				\multirow{3}{*}{\begin{tabular}[c]{@{}c@{}}n=20K, d=20K,\\ r=1K\end{tabular}} & 4.8          & 7.4e-12        & 1230.1     & 2.7e-03     & 408.9      & 7.3e-10    & 914.8      & 2.5e-07      & 376.3  & x     & 506.2                     \\
				& 6.0          & 6.0e-12        & 1833.8     & 2.3e-03     & 1013.2     & 2.1e-09    & 954.2      & 5.4e-07      & 459.6  & x     & 810.2                     \\
				& 10.5         & 1.9e-09        & 1205.7     & 2.2e-03     & 861.3      & 3.0e-09    & 626.5      & 6.9e-07      & 440.9  & x     & 1022.8                    \\ \hline
				\multirow{3}{*}{\begin{tabular}[c]{@{}c@{}}n=40K, d=10K,\\ r=1K\end{tabular}} & 5.0          & 3.4e-12        & 1819.1     & 4.7e-04     & 432.7      & 1.0e-09    & 920.4      & 4.6e-08      & 244.6  & x     & 609.6                     \\
				& 5.9          & 3.7e-12        & 1708.3     & 4.7e-04     & 576.7      & 2.1e-09    & 935.2      & 8.4e-08      & 295.0  & x     & 967.4                     \\
				& 8.5          & 9.0e-10        & 1595.3     & 2.6e-04     & 913.0      & 2.9e-09    & 902.6      & 1.2e-07      & 337.6  & x     & 1183.2                    \\ \hline
				\multirow{3}{*}{\begin{tabular}[c]{@{}c@{}}n=10K, d=40K,\\ r=1K\end{tabular}} & 4.5          & 9.0e-12        & 1825.1     & 1.4e-02     & 332.9      & 6.4e-10    & 932.5      & 1.2e-06      & 503.1  & x     & 388.6                     \\
				& 4.7          & 4.3e-11        & 2648.4     & 1.1e-02     & 887.9      & 1.9e-09    & 941.0      & 2.8e-06      & 592.6  & x     & 652.2                     \\
				& 12.5         & 3.0e-09        & 2327.5     & 1.2e-02     & 821.1      & 3.0e-09    & 616.1      & 3.9e-06      & 629.8  & x     & 829.2                     \\ \hline
			\end{tabular}
		}
		
		\label{table: linear regression compare with benchmarks}
	\end{table}

	\subsection{$\ell_1$-constrained logistic regression}\label{subsection: exp: l1 constrained logistic regression}
	Consider the problem 
	\begin{equation}
		\min_{\bbx\in \R^d} ~~ f(\bbx):=\sum_{i=1}^n \log \Big( 1+ \exp(-b_i \bba_i^\T \bbx) \Big) ~~~ \text{s.t.} ~~~ \| \bbx\|_1 \le C \nonumber
	\end{equation}
	where feature-vector $\bba_i\in \R^d$, response $b_i \in \{-1,1\}$ for $i\in [n]$, and regularization parameter $C>0$. 
	The data is generated as per the underlying model:
	\begin{equation}\label{data generation: logistic}
		\Pb(b_i = 1) = 1/ (1+ \exp(-s\bba_i^\T \bbx^*)),
	\end{equation}
	where, $\bba_1$, \ldots, $\bba_n$ and $\bbx^*$ are generated as in Section~\ref{subsection: exp: l1 constrained least squares}.
	The parameter $s>0$ is used to control the signal-to-noise ratio of the model, and we take $s = 1$.

	We compare PolyCDwA with {FISTA} and AFW. 
	The optimality gap of an algorithm is defined similarly as in Section \ref{subsubsection: linear regression Comparison with benchmarks}. 
	PolyCDwA is run for a maximum of $100$ outer iterations and is terminated earlier if the relative improvement across two successive outer iterations is less than $10^{-9}$.
	AFW and FISTA are run for a maximum of $5000$ outer iterations and are terminated earlier if the relative improvement in the past $50$ iterations is less than $10^{-9}$.

	Table~\ref{table: logistic regression compare with benchmarks} reports the runtimes and optimality gaps (abbreviated as ``gap") of the three algorithms for different values of $n$, $d$ and $r$. The reported results are the average of $5$ independent experiments. For FISTA, only runtimes are reported, as the optimality gaps of FISTA in these examples are all $0$ (i.e. the smallest among the three methods).
	From Table~\ref{table: logistic regression compare with benchmarks}, it can be seen that under our termination rules, FISTA finds the best objective value across all instances. The solutions of PolyCDwA also have high accuracy with optimality gaps $10^{-8} \sim 10^{-10}$, and the runtimes of PolyCDwA are much smaller than FISTA---we improve over FISTA by a factor of $5$X$\sim 10$X. 
	The runtimes of AFW are much longer than those of the other two algorithms. 
	
	\begin{table}[h]
		\centering
		\caption{Comparison of PolyCDwA and benchmarks on $\ell_1$-constrained logistic regression. The terms ``gap'' and ''Nonzeros'' are defined in Table~\ref{table: linear regression compare with benchmarks}.}
		\scalebox{0.71}{
			\begin{tabular}{|c|cc|cc|c|c|} 
				\hline
				\multirow{2}{*}{}                                                            & \multicolumn{2}{c|}{PolyCDwA} & \multicolumn{2}{c|}{AFW} & FISTA & \multirow{2}{*}{Nonzeros}  \\
				& time  & gap                   & time    & gap            & time  &                            \\ 
				\hline
				\multirow{3}{*}{\begin{tabular}[c]{@{}c@{}}n=20K, d=20K\\ r=1K\end{tabular}} & 54.8  & 7.7e-10               & 24635.1 & 6.0e-09        & 657.6 & 624.0                      \\
				& 55.4  & 4.6e-10               & 30918.1 & 5.5e-06        & 695.3 & 1066.0                     \\
				& 111.0 & 7.8e-09               & 27330.0 & 1.0e-03        & 687.9 & 1596.0                     \\ 
				\hline
				\multirow{3}{*}{\begin{tabular}[c]{@{}c@{}}n=40K, d=10K\\ r=1K\end{tabular}} & 51.4  & 6.9e-10               & 22585.3 & 5.3e-09        & 741.8 & 515.0                      \\
				& 45.3  & 5.9e-10               & 32636.4 & 5.1e-06        & 622.1 & 912.0                      \\
				& 98.6  & 1.3e-08               & 27031.4 & 1.3e-03        & 609.5 & 1445.0                     \\ 
				\hline
				\multirow{3}{*}{\begin{tabular}[c]{@{}c@{}}n=10K, d=40K\\ r=1K\end{tabular}} & 85.4  & 9.2e-10               & 21715.4 & 5.7e-09        & 728.6 & 438.0                      \\
				& 73.4  & 3.7e-10               & 32556.7 & 3.4e-06        & 602.8 & 736.0                      \\
				& 113.0 & 1.9e-08               & 26683.0 & 1.6e-03        & 539.7 & 1092.0                     \\
				\hline
			\end{tabular}
		}
		\label{table: logistic regression compare with benchmarks}
	\end{table}

	\subsection{Kernel density estimation}\label{subsection: exp: kernal density estimation}
	We consider the robust kernel density estimation problem \cite{kim2012robust}. 
	Given a set of iid observations $X_1, X_2, \ldots, X_n\in \R^d$ from an underlying distribution with density $g^*(x)$, the goal is to estimate $g^*(\cdot)$ based on $\{X_i\}_1^n$. In particular, let $\cK(\cdot, \cdot): \R^d \times \R^d \rightarrow \R_+$ be a positive semidefinite kernel and let $\cH$ be the Reproducing Kernel Hilbert Space (RKHS) induced by $\cK$. 
	The kernel density estimator of $g^*$ is defined as 
	\begin{equation}\label{problem: KDE1}
		\hat g \in \argmin_{g\in \cH} ~ \sum_{i=1}^n \varphi\lt(\|\cK(\cdot, X_i) - g(\cdot)\|_{\cH}\rt)
	\end{equation}
	where $\varphi: \R \rightarrow \R_+$ is a given loss function, for example, the squared loss $\varphi(t) = t^2$, or the Huber loss 
	\begin{equation}\label{Huber-loss}
		\varphi (t) = \lt\{
		\begin{array}{ll}
			t^2 /2 & ~~\text{if}~~0 \le t \le \mu , \\
			\mu t - \mu^2/2 & ~~\text{otherwise}.
		\end{array}
		\rt.
	\end{equation}
	For the above loss function, by representer theorem \cite{kim2012robust}, there exists $w\in \Dt_n$ such that $\hat g(\cdot) = \sum_{i=1}^n w_i \cK(\cdot, X_i)$. As a result, problem \eqref{problem: KDE1} can be written as 
	\begin{equation}\label{problem: KDE2}
		\small
		\min_{\bbw\in \Dt_n}~ f(\bbw):=	\sum_{i=1}^n \varphi \big(  (  \bbw^\T K \bbw - 2\bbe_i^\T K\bbw + K_{ii}  )^{1/2}  \big) 
	\end{equation}
	where $K \in \R^{n\times n}$ with $K_{ij} = \cK(X_i, X_j)$. 
	We consider problem \eqref{problem: KDE2} with the Huber loss~\eqref{Huber-loss} and Gaussian kernels with a fixed variance $\sigma>0$:
	\begin{equation}\label{gaussian-kernal}
		\small
		\cK_\sigma (\bbx, \bbx') := (2\pi \sigma^2)^{-d/2} \exp \Big(   - {\|\bbx-\bbx'\|^2}/({2\sigma^2} ) \Big) .
	\end{equation} 
	See Appendix~\ref{app-subsection: exp: kernal density estimation} for details on the data generation.

	We compare PolyCDwA with FISTA and AFW on a large example with $d=2$ and $n = 50000$.
	Note that for this example, the matrix $K$ is too large to be maintained in memory, but any entry of $K$ can be computed easily with the formula \eqref{gaussian-kernal}. 
	PolyCDwA is run with a maximum of $30$ (outer) iterations and is terminated earlier if the relative improvement in one iteration is less than $10^{-8}$. 
	FISTA and AFW are run with a maximum of $100$ iterations. 
	We let $f^*$ be the best objective value obtained across all these three methods upon termination, and define the (relative) optimality gap of an algorithm in a manner similar to Section~\ref{subsection: exp: l1 constrained least squares}.

	\begin{figure}[h]
		
		\centering 
		\includegraphics[width=0.5\textwidth]{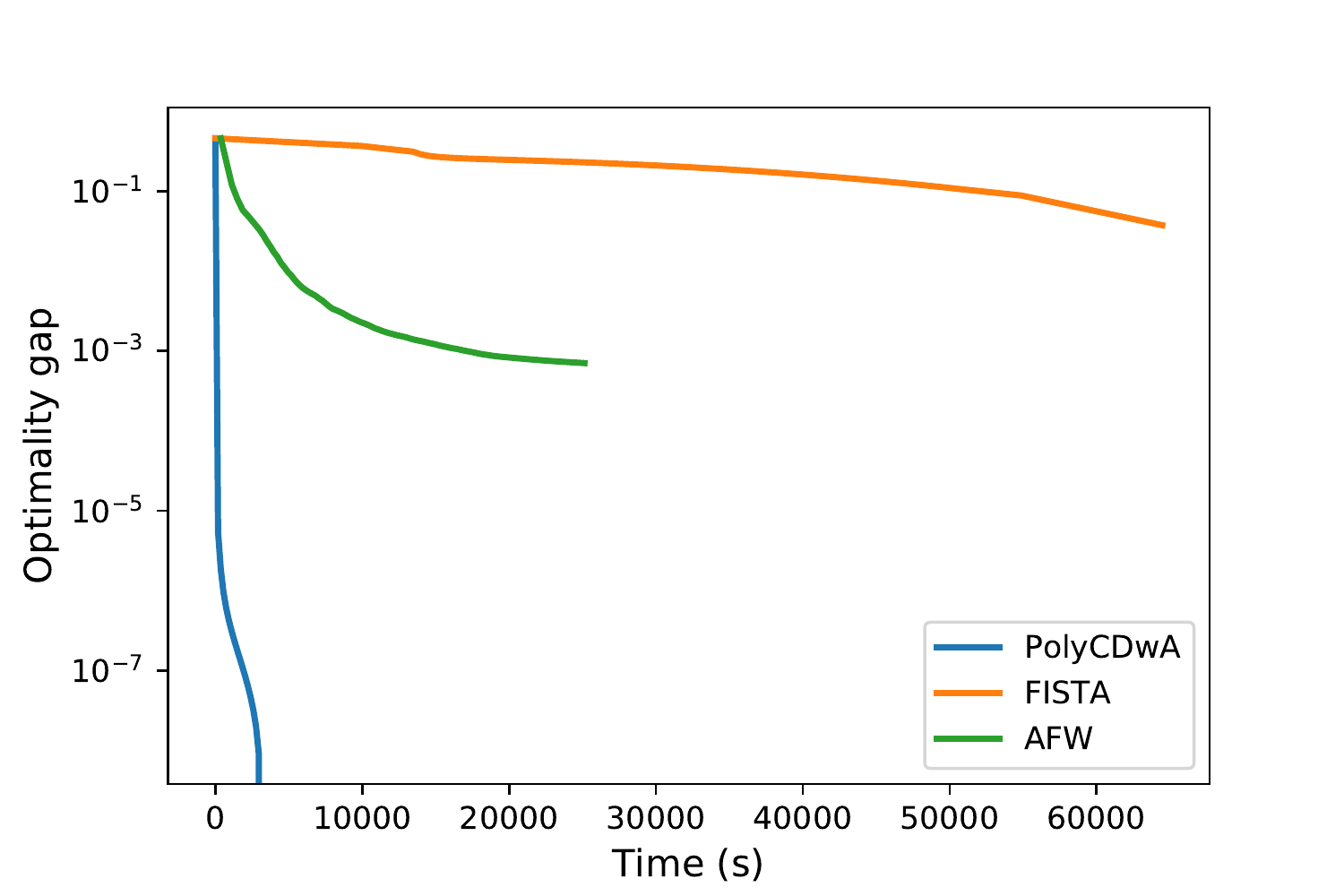}
		
		\caption{Optimality gap vs runtime(s) plotting for PolyCDwA, FISTA and AFW.}
		\label{fig: density compare}
	\end{figure}

	Figure \ref{fig: density compare} presents the performances of PolyCDwA, FISTA, and AFW. For this example, PolyCDwA finds a better objective value upon termination, and is much faster than the others. AFW also makes progress and finds a solution with an optimality gap of approx. $10^{-3}$ in $100$ iterations, but takes much longer time than PolyCDwA. FISTA seems to make slow progress and only finds a solution with optimality gap $10^{-1}\sim 10^{-2}$.

	\appendix

	\section{Proof of Theorem \ref{theorem: PolyCD}}	
	Below we present the proof of Theorem~\ref{theorem: PolyCD}. We first introduce a set of inequalities (which we use in the proof) that follows immediately from the step-size rules.
	
	\begin{lemma}\label{lemma: one-step-results}
		Under the setup of Theorem \ref{theorem: PolyCD}, it holds:
		
		(1) If line search steps \eqref{stepsize: line-search-rule} are used, then for all $t\ge0 $, $i\in [M]$ and all 
		$\bbu\in [ \bbx^{t,i-1},  \bbv^i]$, it holds $\la \na f( \bbx^{t,i}), \bbx^{t,i} - \bbu \ra \le  0 $. 
		
		(2) If 1D gradient steps \eqref{stepsize: 1d-gradient-rule} are used, then for all $t\ge0 $ and $i\in [M]$, we have:
		
		~~~(2.a) 
		For all $u\in [\bbx^{t,i-1},  \bbv^{i}]$, it holds $\la \na f(\bbx^{t,i-1}), \bbx^{t,i} - \bbu \ra \le  L D \| \bbx^{t,i} - \bbx^{t,i-1}\|  $. 
		
		~~~(2.b) $	\la \na f(\bbx^{t,i-1}), \bbx^{t,i} - \bbx^{t,i-1} \ra \le -L \| \bbx^{t,i} - \bbx^{t,i-1}\|^2  $. 
			\end{lemma}
		
		\noindent
		\textit{Proof of Lemma~\ref{lemma: one-step-results}.}
			The proof of (1) follows immediately from the optimality condition of \eqref{stepsize: line-search-rule}. 
			
			To prove (2), 
			by the optimality condition of update~\eqref{stepsize: 1d-gradient-rule}, we have 
			\begin{equation}\label{BCD-opt-condition}
				\la \na f(\bbx^{t,i-1}) + L(\bbx^{t,i} - \bbx^{t,i-1}) , \bbu -\bbx^{t,i}\ra \ge 0 
			\end{equation}
			for all $\bbu\in [\bbx^{t,i-1}, \bbv^{i}] $. Hence 
			\begin{equation}
				\la \na f(\bbx^{t,i-1}), \bbx^{t,i} - \bbu \ra  \le 
				L \la \bbx^{t,i} - \bbx^{t,i-1}, \bbu - \bbx^{t,i} \ra \le LD \| \bbx^{t,i} - \bbx^{t,i-1}\| , 
				\nonumber
			\end{equation}
			this proves $(2.a)$. Taking $\bbu = \bbx^{t,i-1}$ in \eqref{BCD-opt-condition} we have proved $(2.b)$. 
\endproof

	\subsection{Proof of Theorem \ref{theorem: PolyCD} (1)}
	By the convexity of $f(\cdot)$ we know 
	\begin{equation}\label{ineq1-BCD-1}
		f(\bbx^{t,M}) - f(\bbx^*) \le \la \na f(\bbx^{t,M}), \bbx^{t,M} - \bbx^* \ra .
	\end{equation}
	Since $S = \conv(\{\bbv^1,\dots,\bbv^M\})$, there exists $i_* \in [M]$  such that 
	\begin{equation}\label{ineq2-BCD-1}
		\la \na f(\bbx^{t,M}), \bbx^{t,M} - \bbv^{i_*} \ra = \max_{\bbx\in S} ~ \la \na f(\bbx^{t,M}), \bbx^{t,M} - \bbx \ra .
	\end{equation}
	Combining \eqref{ineq1-BCD-1} and \eqref{ineq2-BCD-1} we have 
	\begin{equation}\label{ineq3-BCD-1}
		f(\bbx^{t,M}) - f(\bbx^*) \le 	\la \na f(\bbx^{t,M}), \bbx^{t,M} - \bbv^{i_*} \ra. 
	\end{equation}
	
	By using equality \eqref{shortcut-eq1} in Lemma~\ref{lemma: reduction} with $i = i_*$, $j=M$ and $\bbz = \bbv^{i_*}$, we have 
	\begin{equation}\label{ineq4-BCD-1}
		\begin{aligned}
			\la \na f(\bbx^{t,M}), \bbx^{t,M} - \bbv^{i_*} \ra  = \la \na f(\bbx^{t,i_*}), \bbx^{t,i_*} - \bbv^{i_*} \ra 
			+ \sum_{k=i_*+1}^M \la \na f(\bbx^{t,k}), \bbx^{t,k} - \bbx^{t,k-1} \ra  \\
			+ \sum_{k=i_*+1}^M \la \na f(\bbx^{t,k}) - \na f(\bbx^{t,k-1}) , \bbx^{t,k-1} - \bbv^{i_*} \ra .
		\end{aligned}
	\end{equation}
	By Lemma \ref{lemma: one-step-results} (1) and noting that 
	$\bbx^{t,k-1} \in [\bbx^{t,k-1}, \bbv^{k}]$ and $ \bbv^{i_*} \in [ \bbx^{t,i_*-1}, \bbv^{i_*}] $, 
	we have 
	\begin{equation}\label{ineq5-BCD-1}
		\la \na f(\bbx^{t,i_*}), \bbx^{t,i_*} - \bbv^{i_*} \ra \le 0, ~~\text{and}~~~ \la \na f( \bbx^{t,k}), \bbx^{t,k} - \bbx^{t,k-1} \ra \le 0 ~~ \forall k\in [M]
	\end{equation}
	which shows that the first two terms in the rhs of~\eqref{ineq4-BCD-1} are nonpositive. Using this observation and combining \eqref{ineq3-BCD-1}, \eqref{ineq4-BCD-1}, we get:
	\begin{equation}
		\begin{aligned}
			f( \bbx^{t,M}) - f( \bbx^*) \le \ & \sum_{k=i_*+1}^M \la \na f( \bbx^{t,k}) - \na f( \bbx^{t,k-1}) , \bbx^{t,k-1} - \bbv^{i_*} \ra \\
			\le \ &
			\sum_{k=i_*+1}^M   D \| \na f  ( \bbx^{t,k}) - \na f( \bbx^{t,k-1} )\|  
		\end{aligned}
		\nonumber
	\end{equation}
	where the second inequality is because $\|\bbx^{t,k-1} - \bbv^{i_*} \| \le D$ (by the definition of $D$ and $\bbx^{t,k-1}, \bbv^{i_*} \in S$). Squaring both sides of the inequality above, 
	\begin{equation}\label{ineq7-BCD-1}
		\begin{aligned}
			(f(\bbx^{t,M}) - f(\bbx^*) )^2 \le \ & D^2 \Big( \sum_{k=1}^M \| \na f  (\bbx^{t,k}) - \na f( \bbx^{t,k-1} )\|  \Big)^2 \\
			\le \ & 
			MD^2  \sum_{k=1}^M \| \na f  (\bbx^{t,k}) - \na f( \bbx^{t,k-1} )\|^2 
		\end{aligned}
	\end{equation}
	where the second inequality is by Jensen's inequality. 
	On the other hand, using Lemma~\ref{lemma: convex-basic} 
	we have 
	\begin{equation}
		\begin{aligned}
			f(\bbx^{t,k-1}) - f(\bbx^{t,k}) \ge  \ & \la \na f(\bbx^{t,k}), \bbx^{t,k-1} - \bbx^{t,k} \ra + (1/(2L))
			\| \na f  (\bbx^{t,k}) - \na f( \bbx^{t,k-1} )\|^2  \\
			\ge \ & 
			(1/(2L))
			\| \na f  (\bbx^{t,k}) - \na f( \bbx^{t,k-1} )\|^2 
		\end{aligned}
		\nonumber
	\end{equation}
	where the second inequality makes use of \eqref{ineq5-BCD-1}. Summing up the above inequality across $k\in [M]$, one has 
	\begin{equation}\label{ineq9-BCD-1}
		2L(f(\bbx^{t,0}) - f(\bbx^{t,M} ) )\ge \sum_{k=1}^M \| \na f  (\bbx^{t,k}) - \na f( \bbx^{t,k-1} )\|^2 .
	\end{equation}
	Combining \eqref{ineq7-BCD-1} and \eqref{ineq9-BCD-1}, we have 
	\begin{equation}
		(f(\bbx^{t,M}) - f(\bbx^*) )^2  \le 2MLD^2 (f(\bbx^{t,0}) - f(\bbx^{t,M})). \nonumber
	\end{equation}
	Recall that $\bbx^{t} = \bbx^{t,0}$ and $\bbx^{t+1} = \bbx^{t,M}$, so we have 
	\begin{equation}
		(f(\bbx^{t+1}) - f(\bbx^*) )^2  \le 2MLD^2 (f(\bbx^{t}) - f(\bbx^{t+1})). \nonumber
	\end{equation}
	Using Lemma \ref{lemma: technical1} with $a_t = f(\bbx^t) - f(\bbx^*)$, the proof is complete.

	\subsection{Proof of Theorem \ref{theorem: PolyCD} (2)}
	By the same argument from \eqref{ineq1-BCD-1} to \eqref{ineq3-BCD-1}, we have 
	\begin{equation}\label{ineq3-BCD-1-rep}
		f(\bbx^{t,M}) - f(\bbx^*) \le 	\la \na f(\bbx^{t,M}), \bbx^{t,M} - \bbv^{i_*} \ra. 
	\end{equation}
	Using equation \eqref{shortcut-eq2} in Lemma~\ref{lemma: reduction} with $i = i_*$, $j=M$ and $\bbz = \bbv^{i_*}$, we have 
	\begin{equation}\label{ineq4-BCD}
		\begin{aligned}
			&
			\la \na f(\bbx^{t,M}), \bbx^{t,M} - \bbv^{i_*} \ra \\
			= \ &
			\la \na f(\bbx^{t,i_*-1}), \bbx^{t,i_*} - \bbv^{i_*} \ra  + \sum_{k=i_*+1}^M \la \na f(\bbx^{t,k-1}), \bbx^{t,k} - \bbx^{t,k-1} \ra \\
			& + \sum_{k=i_*}^M \la \na f(\bbx^{t,k}) - \na f(\bbx^{t,k-1}) , \bbx^{t,k} - \bbv^{i_*} \ra .
		\end{aligned}
	\end{equation}
	By Lemma \ref{lemma: one-step-results} (2.a) and noting that $ \bbv^{i_*}  \in [\bbx^{t,i_*-1},  \bbv^{i_*}] $, 
	we have 
	\begin{equation}\label{ineq5-BCD}
		\la \na f(\bbx^{t,i_*-1}), \bbx^{t,i_*} - \bbv^{i_*} \ra \le L D \| \bbx^{t,i_*} - \bbx^{t,i_*-1} \| .
	\end{equation}
	By Lemma \ref{lemma: one-step-results} (2.b) 
	we have 
	\begin{equation}\label{ineq6-BCD}
		\la \na f( \bbx^{t,k-1}), \bbx^{t,k} - \bbx^{t,k-1} \ra \le -L \| \bbx^{t,k} - \bbx^{t,k-1}\|^2  \le 0 .
	\end{equation}
	Combining \eqref{ineq3-BCD-1-rep}, \eqref{ineq4-BCD}, \eqref{ineq5-BCD} and \eqref{ineq6-BCD} we have 
	\begin{equation}\label{ineq6pt5-BCD}
		\begin{aligned}
			f(\bbx^{t,M}) - f(\bbx^*) \le \ &  LD \| \bbx^{t,i_*} -\bbx^{t,i_*-1} \| + \sum_{k=i_*}^M \la \na f(\bbx^{t,k}) - \na f(\bbx^{t,k-1}) , \bbx^{t,k} - \bbv^{i_*} \ra \\
			\le \ &
			LD \| \bbx^{t,i_*} -\bbx^{t,i_*-1} \| +  \sum_{k=i_*}^M  L D \| \bbx^{t,k} - \bbx^{t,k-1} \|  
			\\
			\le \ &
			2LD \sum_{k=1}^M \| \bbx^{t,k} - \bbx^{t,k-1} \|  
		\end{aligned}
	\end{equation}
	where the second inequality makes use of the $L$-smoothness of $f(\cdot)$. 
	Squaring the two extreme sides of~\eqref{ineq6pt5-BCD}, and using Jensen's inequality, we have:
	\begin{equation}\label{ineq7-BCD}
		(f(\bbx^{t,M}) - f(\bbx^*) )^2 \le 
		4ML^2D^2  \sum_{k=1}^M \| \bbx^{t,k} - \bbx^{t,k-1} \|^2 . 
	\end{equation}
	On the other hand, by \eqref{ineq: convex-basic1} in Lemma~\ref{lemma: convex-basic} we have 
	\begin{equation}
		\begin{aligned}
			f(\bbx^{t,k-1}) - f(\bbx^{t,k}) \ge \ & \la \na f(\bbx^{t,k-1}), \bbx^{t,k-1} - \bbx^{t,k} \ra - (L/2)\| \bbx^{t,k-1} - \bbx^{t,k} \|^2 \\
			\ge \ &
			(L/2)\| \bbx^{t,k-1} - \bbx^{t,k} \|^2 
		\end{aligned}
		\nonumber
	\end{equation}
	where the second inequality is by Lemma \ref{lemma: one-step-results} (2.b). Summing the above inequality over $k\in [M]$, we have 
	\begin{equation}\label{ineq8-BCD}
		f(\bbx^{t,0}) - f(\bbx^{t,M}) \ge (L/2)  \sum_{k=1}^M \| \bbx^{t,k} - \bbx^{t,k-1} \|^2  .
	\end{equation}
	Combining \eqref{ineq7-BCD} and \eqref{ineq8-BCD}, we have 
	\begin{equation}
		(f(\bbx^{t,M}) - f(\bbx^*) )^2  \le 8MLD^2 (f(\bbx^{t,0}) - f(\bbx^{t,M})). \nonumber
	\end{equation}
	Recall that $\bbx^{t} = \bbx^{t,0}$ and $\bbx^{t+1} = \bbx^{t,M}$, so we have 
	\begin{equation}
		(f(\bbx^{t+1}) - f(\bbx^*) )^2  \le 8MLD^2 (f(\bbx^{t}) - f(\bbx^{t+1})). \nonumber
	\end{equation}
	Using Lemma \ref{lemma: technical1} with $a_t = f(\bbx^t) - f(\bbx^*)$, the proof is complete.

	\section{Proof of Theorem \ref{theorem: PolyCDwA strong convex}}\label{section: proof of theorem PolyCDwA strong convex}
	
	Let $
	V^{t,i} := \{   j \in [M] ~|~ \lam^{t,i}_j >0    \} $, and we call $V^{t,i}$ the \textit{vertex-support} of the iterate $\bbx^{t,i}$. 
	We first prove a few technical lemmas that are used in proofs of both parts (1) and (2) and then proceed with proving these two parts.

	\begin{lemma}\label{lemma: key bound away}
		Suppose the assumptions in the statement of Theorem \ref{theorem: PolyCDwA strong convex} hold. 
		Let $\{\bbx^{t,i}\} $ and $\{\bblam^{t,i}\}$ ($t\ge 0, 0\le  i  \le M$) be the sequences generated by Algorithm~\ref{alg: PolyCDwA} with either line-search steps  \eqref{stepsize: line-search-rule-away} or 1D gradient steps  \eqref{stepsize: 1D-gradient-rule-away}.
		Then there exists $\eta\in [0,2]$, and indices $j_1 \in V^{t,M}$ and $j_2 \in [M]$ 
		such that 
		\begin{equation}
			f(\bbx^{t,M}) - f(\bbx^*) 
			\le 
			\frac{\eta}{2}    \la \na f(\bbx^{t,M}), \bbv^{j_1} - \bbv^{j_2} \ra  - \frac{\mu \psi_S^2}{8} \eta^2 . \nonumber
		\end{equation}
	\end{lemma}
\noindent
\textit{Proof of Lemma~\ref{lemma: key bound away}.}
		Since $f(\cdot)$ is $\mu$-strongly convex, 
		we have 
		\begin{equation}\label{g1}
			f(\bbx^{t,M}) - f(\bbx^*) \le  \la \na f(\bbx^{t,M}) , \bbx^{t,M} - \bbx^* \ra - \frac{\mu}{2} \| \bbx^{t,M} - \bbx^* \|^2 .
		\end{equation}
		Using Lemma \ref{lemma: poly} with $\bbx = \bbx^{t,M}$ and $\bby =  \bbx^*$, there exists $\bblam^* \in \Dt_M$ such that $\bbx^* = \sum_{i=1}^M \lam^*_i \bbv^i $,  and 
		\begin{equation}\label{norm-bound}
			\| \bbx^{t,M} - \bbx^* \|^2 \ge (\psi_S/2)^2 \| \bblam^{t,M} - \bblam^* \|_1^2 .
		\end{equation}
		Let $\eta := \| \bblam^{t,M} - \bblam^* \|_1 $, then we have $\eta \le \| \bblam^{t,M} \|_1+\| \bblam^* \|_1 \le 2$. By Lemma \ref{lemma: simplex vector decompose}, there exist $\wtd \bbp, \wtd \bbq \in \Dt_M$ such that $\supp(\wtd \bbp) \subseteq V^{t,M}$ and 
		\begin{equation}\label{decompose1}
			\bblam^{t,M} - \bblam^* = \frac{\eta}{2} (\wtd \bbp - \wtd \bbq). 
		\end{equation}
		Define $A = [\bbv^1,\dots,\bbv^M] \in \R^{d\times M}$. Then we have 
		\begin{equation}\label{formula-m1}
			\bbx^{t,M} - \bbx^* = A(\bblam^{t,M} - \bblam^*)  = 
			\frac{\eta}{2} A (\wtd \bbp - \wtd \bbq)  
		\end{equation}
		where the second equality is because of \eqref{decompose1}. 
		Recall that $\wtd \bbp, \wtd \bbq \in \Dt_M$ and $\supp(\wtd \bbp) \subseteq V^{t,M}$, and {because a linear function achieves its maximal and minimal values at extreme points of a polytope}, we know that
		there exist $j_1 \in V^{t,M}$ and $j_2 \in [M]$ such that 
		\begin{equation}\label{push-to-ext1}
			\la \na f(\bbx^{t,M}) , A \wtd \bbp \ra \le 
			\la \na f(\bbx^{t,M}) , \bbv^{j_1} \ra, \quad 
			\la \na f(\bbx^{t,M}) , A \wtd \bbq \ra \ge 
			\la \na f(\bbx^{t,M}) , \bbv^{j_2} \ra.
		\end{equation}
		Combining the two inequalities in \eqref{push-to-ext1} and using
		\eqref{formula-m1}, we have
		\begin{equation}\label{formula-m2}
			\la   \na f(\bbx^{t,M}),   \bbx^{t,M} - \bbx^*  \ra \le (\eta/2) 	\la   \na f(\bbx^{t,M}), \bbv^{j_1} - \bbv^{j_2} \ra .
		\end{equation}
		Combining \eqref{g1}, \eqref{norm-bound}, \eqref{formula-m2} and recalling that $\eta = \| \bblam^{t,M} - \bblam^* \|_1$, the proof is complete. 
\endproof

	\subsection{Proof of Theorem \ref{theorem: PolyCDwA strong convex} (1)}
	We first present a few technical lemmas. 
	\begin{lemma}\label{lemma: line-search-results-away}
		Suppose line-search steps \eqref{stepsize: line-search-rule-away} are used, 
		then for all $t\ge 0$ and $i\in [M]$, 
		(a) $\la \na f(\bbx^{t,i}) , \bbx^{t,i} - \bbx^{t,i-1} \ra \le 0$. 
		(b) $\la \na f(\bbx^{t,i}) , \bbx^{t,i} - \bbv^i \ra  \le 0$.
		
	\end{lemma}
	The correctness of Lemma~\ref{lemma: line-search-results-away} can be immediately verified by the optimality condition of the line-search steps
	\eqref{stepsize: line-search-rule-away}. 
	
	For any given $t\ge0 $ and integers $j_1 \in V^{t,M}$ and $j_2 \in [M]$, we define the following quantities:
	\begin{equation}\label{def:J1-J2}
		J_1 := \la \na f(\bbx^{t,j_1}), \bbv^{j_1} - \bbx^{t,j_1} \ra ,\quad 
		J_2 := \la \na f(\bbx^{t,j_1}), \bbx^{t,j_1} - \bbx^{t,j_2} \ra
	\end{equation}
	\begin{equation}\label{def:J3}
		J_3 := \la \na f(\bbx^{t,j_2}), \bbx^{t,j_2} - \bbv^{j_2} \ra
	\end{equation}
	\begin{equation}\label{def:J4}
		J_4 := \la \na f(\bbx^{t,M}) - \na f(\bbx^{t,j_1}), \bbv^{j_1} - \bbx^{t,j_2} \ra + 
		\la \na f(\bbx^{t,M}) - \na f(\bbx^{t,j_2}), \bbx^{t,j_2} - \bbv^{j_2} \ra
	\end{equation}
	Note that $J_1 $--$ J_4$ depend on $j_1,j_2$ and $t$; For notational simplicity we drop the dependence on $j_1,j_2$ and $t$. 
	We have the following lemma for upper bounding $J_1$--$J_4$. 
	\begin{lemma}\label{lemma: bound J1-J5}
		For any $t\ge 0$, $j_1 \in V^{t,M}$ and $j_2 \in [M]$ with $j_1\neq j_2$, the quantities defined in \eqref{def:J1-J2}, \eqref{def:J3} and \eqref{def:J4} satisfy:
		
		(1) $J_1 =0.$ 
		
		(2) $J_2 \le \sum_{k=1}^M  \la  \na f(\bbx^{t,k}), \bbx^{t,k-1} - \bbx^{t,k} \ra +
		\sum_{k=1}^M D\| \na f(\bbx^{t,k-1})   - \na f(\bbx^{t,k})  \| $. 
		
		(3) $J_3 \le 0 $. 
		\quad
		(4) $J_4 \le  2D \sum_{k=1}^M \| \na f(\bbx^{t,k-1})   - \na f(\bbx^{t,k})  \|  $. 
	\end{lemma}

		\noindent
		\textit{Proof of Lemma~\ref{lemma: bound J1-J5}}
			(1) We discuss different cases of $\al_{t,i}$. $(i)$ If $\al_{t,j_1} = 1 $, then $\bbx^{t,j_1} = \bbv^{j_1}$, and hence $ J_1= 0$. 
			$(ii)$ If $\al_{t,j_1}\in ( -\ga_{t, j_1} , 1)$, then by \eqref{stepsize: line-search-rule-away}, we know that $ \na f(\bbx^{t,j_1}) $ should be orthogonal to the vector $\bbv^{j_1} - \bbx^{t,j_1} $, i.e., $ \la \na f(\bbx^{t,j_1}) , \bbv^{j_1} - \bbx^{t,j_1} \ra = 0$. 
			$(iii)$ We consider the case when $\al_{t,j_1} = -\ga_{t, j_1}$, and show by contradiction that this cannot happen. By the definition of $\ga_{t, j_1}$ and the first equality in \eqref{PolyCDwA update lam}, we know
			$ \lam_{j_1}^{t,j_1} = 0 $. As a result, by the second equality in \eqref{PolyCDwA update lam} we know $ \lam_{j_1}^{t,j_1+1} = \lam_{j_1}^{t,j_1+2} = \cdots = \lam_{j_1}^{t,M} = 0$. But this is contradictory to the fact that $ j_1 \in V^{t,M} = \{   j \in [M] ~|~ \lam^{t,M}_j >0    \}$, so the case $\al_{t,j_1} = -\ga_{t, j_1} $ cannot arise. 
			
			(2) To bound $J_2$, we consider 2 cases. 
			$(i)$ If $j_1>j_2$, then 
			using \eqref{shortcut-eq1} in Lemma~\ref{lemma: reduction} with $\bbz = \bbx^{t,j_2}$, $i = j_2$ and $j=  j_1$ we have 
			\begin{equation}
				\begin{aligned}
					J_2 = \ &  \sum_{k=j_1+1}^{j_2} \lt(  \la  \na f(\bbx^{t,k}), \bbx^{t,k} - \bbx^{t,k-1} \ra + \la \na f(\bbx^{t,k})   - \na f(\bbx^{t,k-1})    ,   \bbx^{t,k-1} -\bbx^{t,j_2} \ra \rt) \\
					\mathop{\le} \ &
					\sum_{k=j_1+1}^{j_2}   \la \na f(\bbx^{t,k})   - \na f(\bbx^{t,k-1})    ,   \bbx^{t,k-1} -\bbx^{t,j_2} \ra \\
					\le \ & 
					\sum_{k=j_1+1}^{j_2}   D \| \na f(\bbx^{t,k})   - \na f(\bbx^{t,k-1})      \| \le 
					\sum_{k=1}^{M}   D \| \na f(\bbx^{t,k})   - \na f(\bbx^{t,k-1})      \| .
				\end{aligned}
				\nonumber
			\end{equation}
			where the first inequality makes use of Lemma \ref{lemma: line-search-results-away} (a); the second inequality makes use of $\|\bbx^{t,k-1} -\bbx^{t,j_2}\| \le D$.
			$(ii)$ If $j_1 < j_2$, then using \eqref{shortcut-eq1} in Lemma~\ref{lemma: reduction} with $\bbz = \bbx^{t,j_2}$, $i = j_1$ and $j=  j_2$ we have 
			\begin{equation}
				\begin{aligned}
					J_2 = \ & \sum_{k=j_2+1}^{j_1} \lt(  \la  \na f(\bbx^{t,k}), \bbx^{t,k-1} - \bbx^{t,k} \ra + \la \na f(\bbx^{t,k-1})   - \na f(\bbx^{t,k})    ,   \bbx^{t,k-1} - \bbx^{t,j_2} \ra \rt) \\
					\le \ &
					\sum_{k=j_2+1}^{j_1} \lt(  \la  \na f(\bbx^{t,k}), \bbx^{t,k-1} - \bbx^{t,k} \ra + D\| \na f(\bbx^{t,k-1})   - \na f(\bbx^{t,k})  \| \rt) \\
					\le \ &
					\sum_{k=1}^M \lt(  \la  \na f(\bbx^{t,k}), \bbx^{t,k-1} - \bbx^{t,k} \ra + D\| \na f(\bbx^{t,k-1})   - \na f(\bbx^{t,k})  \| \rt) 
				\end{aligned}
				\nonumber
			\end{equation}
			where the last inequality is because $\la  \na f(\bbx^{t,k}), \bbx^{t,k-1} - \bbx^{t,k} \ra \ge 0 $ (by Lemma \ref{lemma: line-search-results-away} (a)). 
			Combining these 2 cases, we have proved (2). 
			
			(3) 
			The conclusion follows immediately from Lemma \ref{lemma: line-search-results-away} (b). 
			
			(4) To  bound $J_4$, we have 
			\begin{equation}
				\begin{aligned}
					J_4  = \ & \la \na f(\bbx^{t,M}) - \na f(\bbx^{t,j_1}), \bbv^{j_1} - \bbx^{t,j_2} \ra + 
					\la \na f(\bbx^{t,M}) - \na f(\bbx^{t,j_2}), \bbx^{t,j_2} - \bbv^{j_2} \ra \\
					\le \ & 
					D \| \na f(\bbx^{t,M}) - \na f(\bbx^{t,j_1})   \|  + D \| \na f(\bbx^{t,M}) - \na f(\bbx^{t,j_2})  \| \\
					\le \ &
					2D \sum_{k=1}^M \| \na f(\bbx^{t,k-1})   - \na f(\bbx^{t,k})  \| .
				\end{aligned}
				\nonumber
			\end{equation}
			This completes the proof of Lemma~\ref{lemma: bound J1-J5}. 
	\endproof

	The results of Lemma~\ref{lemma: bound J1-J5} immediately yields the following lemma. 
	\begin{lemma}\label{lemma:key1}
		For any $t\ge 0$, $j_1 \in V^{t,M}$ and $j_2 \in [M]$ with $j_1 \neq j_2$, it holds 
		\begin{equation}
			\la \na f(\bbx^{t,M}), \bbv^{j_1} - \bbv^{j_2} \ra \le 
			\sum_{k=1}^M  \la  \na f(\bbx^{t,k}), \bbx^{t,k-1} - \bbx^{t,k} \ra + 
			3D \sum_{k=1}^M \| \na f(\bbx^{t,k-1})   - \na f(\bbx^{t,k})  \| .
			\nonumber
		\end{equation}
		\end{lemma}
	\noindent
	\textit{Proof of Lemma~\ref{lemma:key1}.}
			By some algebra, 
			we have the decomposition:
			\begin{equation}
				\begin{aligned}
					&
					\la \na f(\bbx^{t,M}), \bbv^{j_1} - \bbv^{j_2} \ra \\
					= \ &
					\la \na f(\bbx^{t,M}), \bbv^{j_1} - \bbx^{t,j_1} \ra + 
					\la \na f(\bbx^{t,M}), \bbx^{t,j_1} - \bbx^{t,j_2} \ra +
					\la \na f(\bbx^{t,M}), \bbx^{t,j_2} - \bbv^{j_2} \ra  
					\\
					= \ &
					J_1 + J_2 + J_3 + J_4  .
				\end{aligned}
				\nonumber
			\end{equation}
			Using the above equality and Lemma~\ref{lemma: bound J1-J5}, the conclusion is reached.
\endproof

	With Lemmas~\ref{lemma: line-search-results-away} and \ref{lemma:key1} at hand, we are ready to present the proof of Theorem~\ref{theorem: PolyCDwA strong convex} (1). 
	Denote $\wtd \mu := \psi_S^2/4$. Then
	by Lemma \ref{lemma: key bound away}, there exists $\eta\in [0,2]$, and integers $j_1 \in V^{t,M}$ and $j_2 \in [M]$ 
	such that 
	\begin{equation}\label{key3}
		f(\bbx^{t,M}) - f(\bbx^*) 
		\le 
		\frac{\eta}{2}    \la \na f(\bbx^{t,M}), \bbv^{j_1} - \bbv^{j_2} \ra  - \frac{\wtd \mu}{2} \eta^2  .
	\end{equation}
	We can assume $j_1\neq j_2$, since otherwise by \eqref{key3} we have $f(\bbx^{t,M}) \le f(\bbx^*)$, and the conclusion of Theorem~\ref{theorem: PolyCDwA strong convex} (1) holds true trivially. 
	Making use of Lemma~\ref{lemma:key1} and by \eqref{key3}, we have
	\begin{equation}\label{key5}
		\begin{aligned}
			&
			f(\bbx^{t,M}) - f(\bbx^*)  \\
			\le \ & \frac{\eta}{2}  \sum_{k=1}^M  \la  \na f(\bbx^{t,k}), \bbx^{t,k-1} - \bbx^{t,k} \ra  + \frac{3\eta D}{2} \sum_{k=1}^M \| \na f(\bbx^{t,k-1})   - \na f(\bbx^{t,k})  \|  - \frac{\wtd \mu}{2} \eta^2 \\
			\le \ &
			\sum_{k=1}^M  \la  \na f(\bbx^{t,k}), \bbx^{t,k-1} - \bbx^{t,k} \ra  + \frac{3\eta D}{2} \sum_{k=1}^M \| \na f(\bbx^{t,k-1})   - \na f(\bbx^{t,k})  \|  - \frac{\wtd \mu}{2} \eta^2 
		\end{aligned}
	\end{equation}
	where the second inequality is because $\eta \le 2 $. By Cauchy-Schwarz inequality we have 
	\begin{equation}\label{key6}
		\begin{aligned}
			&
			\frac{3\eta D}{2} \sum_{k=1}^M \| \na f(\bbx^{t,k-1})   - \na f(\bbx^{t,k})  \|  - \frac{\wtd \mu}{2} \eta^2 \\
			\le \ & \frac{1}{2\wtd \mu} (3D/2)^2 \Big( \sum_{k=1}^M \| \na f(\bbx^{t,k-1})   - \na f(\bbx^{t,k})  \|   \Big)^2 \\
			\le \ &
			\frac{9MD^2}{8\wtd \mu} \sum_{k=1}^M \| \na f(\bbx^{t,k-1})   - \na f(\bbx^{t,k})  \|^2  
		\end{aligned}
	\end{equation}
	where the second inequality is by Jensen's inequality. 
	Combining \eqref{key5} and \eqref{key6} we have 
	\begin{equation}\label{to-final-1}
		f(\bbx^{t,M}) - f(\bbx^*) \le  \sum_{k=1}^M  \la  \na f(\bbx^{t,k}), \bbx^{t,k-1} - \bbx^{t,k} \ra  + \frac{9MD^2}{8\wtd \mu} \sum_{k=1}^M \| \na f(\bbx^{t,k-1})   - \na f(\bbx^{t,k})  \|^2  .
	\end{equation}
	On the other hand, by the $L$-smoothness of $f(\cdot)$ and inequality \eqref{ineq: convex-basic2}, 
	we have 
	\begin{equation}\label{to-final-1pt05}
		f(\bbx^{t, i-1}) - f(\bbx^{t, i}) \ge \la \na f(\bbx^{t,i}) , \bbx^{t,i-1} - \bbx^{t,i} \ra + \frac{1}{2L} \| \na f(\bbx^{t,i}) - \na f(\bbx^{t,i-1}) \|^2  
	\end{equation}
	for all $i\in [M]$. Summing~\eqref{to-final-1pt05} over $i\in[M]$,  we have 
	\begin{equation}\label{to-final-2}
		\begin{aligned}
			&
			f(\bbx^{t,0}) - f(\bbx^{t,M}) \\
			\ge \ & 
			\sum_{k=1}^M  \la  \na f(\bbx^{t,k}), \bbx^{t,k-1} - \bbx^{t,k} \ra
			+ \frac{1}{2L }     \sum_{k=1}^M \| \na f(\bbx^{t,k-1})   - \na f(\bbx^{t,k})  \|^2  \\
			\ge \ &
			\max\Big\{ \sum_{k=1}^M  \la  \na f(\bbx^{t,k}), \bbx^{t,k-1} - \bbx^{t,k} \ra,  \frac{1}{2L }     \sum_{k=1}^M \| \na f(\bbx^{t,k-1})   - \na f(\bbx^{t,k})  \|^2 \Big\}
		\end{aligned}
	\end{equation}
	By \eqref{to-final-1} and \eqref{to-final-2} we have 
	\begin{equation}
		\begin{aligned}
			f(\bbx^{t,M}) - f(\bbx^*) 
			\le & f(\bbx^{t,0}) - f(\bbx^{t,M}) + \frac{9MLD^2}{4\wtd\mu} (f(\bbx^{t,0}) - f(\bbx^{t,M}) ) \\
			= & G (f(\bbx^{t,0}) - f(\bbx^{t,M}))
		\end{aligned}
		\nonumber
	\end{equation}
	where the last equality is because $1 + 9MLD^2 /(4\wtd \mu) = 1+ 9MLD^2 /(\mu \psi_S^2)  =G $. 
	Recall that $ \bbx^{t+1} = \bbx^{t,M} $ and $\bbx^t = \bbx^{t,0}$, we have 
	\begin{equation}
		f(\bbx^{t+1}) - f^* \le G (f(\bbx^t) - f(\bbx^{t+1}))  \nonumber
	\end{equation}
	or equivalently 
	\begin{equation}
		f(\bbx^{t+1}) - f^* \le \frac{G}{1+G} ( f(\bbx^t) - f^* ),   \nonumber
	\end{equation}
	which completes the proof of Theorem~\ref{theorem: PolyCDwA strong convex} (1).

	\subsection{Proof of Theorem \ref{theorem: PolyCDwA strong convex} (2)}
	We first prove a few technical lemmas. 
	\begin{lemma}\label{lemma: proximal-gradient-results-away}
		Suppose 1D-gradient steps \eqref{stepsize: 1D-gradient-rule-away} are used. 
		Then for all $t\ge 0$ and $ i\in [M]$, 
		
		(a) $\la \na  f(\bbx^{t,i-1}), \bbx^{t,i} - \bbx^{t,i-1} \ra \le - L \| \bbx^{t,i} - \bbx^{t,i-1} \|^2$. 
		
		(b) $\la \na f(\bbx^{t,i-1}), \bbx^{t,i} - \bbv^i  \ra  \le LD \| \bbx^{t,i} - \bbx^{t,i-1} \| $. 
		
		(c) If $\al_{t,i} \in (-\ga_{t,i},1) $, then $\la \na f(\bbx^{t,i-1}), \bbv^i - \bbx^{t,i} \ra \le LD \| \bbx^{t,i} - \bbx^{t,i-1} \| $. 
	\end{lemma}
	\noindent
	\textit{Proof of Lemma~\ref{lemma: proximal-gradient-results-away}.}
		Let $ \cL := \{ \bbx^{t,i-1} + \al ( \bbv^i - \bbx^{t,i-1}) ~|~ \al\in [-\ga_{t,i},1] \}$. (Note that $\cL$ depends on $t,i$; we drop the dependence on $t,i$ for notational convenience). 
		By the optimality condition of \eqref{stepsize: 1D-gradient-rule-away} we have 
		\begin{equation}\label{pgd-opt-condition}
			\la \na  f( \bbx^{t,i-1}) +L( \bbx^{t,i} - \bbx^{t,i-1} ), \bbu - \bbx^{t,i}\ra \ge 0~~~ \forall ~ \bbu\in \cL. 
		\end{equation}
		
		(a) {Letting $ \bbu = \bbx^{t,i-1}$ in~\eqref{pgd-opt-condition}, we obtain the conclusion of 
			part~(a).} 
		
		(b) {Letting $\bbu = \bbv^i$ in~\eqref{pgd-opt-condition}, we have} 
		\begin{equation}
			\la \na f(\bbx^{t,i-1}), \bbx^{t,i} - \bbv^i \ra \le 
			L  \la \bbx^{t,i} - \bbx^{t,i-1}, \bbv^i - \bbx^{t,i} \ra \le 
			LD \| \bbx^{t,i} - \bbx^{t,i-1} \| 
			\nonumber
		\end{equation}
		which completes the proof of part~(b).
		
		(c)	If $\al_{t,i} \in (-\ga_i,1)$ is true, then the inequality in \eqref{pgd-opt-condition} holds as equality. As a result, 
		\begin{equation}
			\la \na f(\bbx^{t,i-1}), \bbv^i - \bbx^{t,i} \ra = 
			L  \la \bbx^{t,i} - \bbx^{t,i-1}, \bbx^{t,i} - \bbv^i\ra \le 
			LD \| \bbx^{t,i} - \bbx^{t,i-1} \|
			\nonumber
		\end{equation}
		which completes the proof of part (c).
\endproof

	For any given $t\ge0 $ and integers $j_1 \in V^{t,M}$ and $j_2 \in [M]$, we define the following quantities:
	\begin{equation}\label{def:J1-J2_v2}
		\td J_1 := \la \na f(\bbx^{t,M}) - \na f(\bbx^{t,j_1-1}), \bbv^{j_1} - \bbx^{t,j_2} \ra ,\quad 
		\td J_2 := \la \na f(\bbx^{t,j_1-1}), \bbx^{t,j_1} - \bbx^{t,j_2} \ra
	\end{equation}
	\begin{equation}\label{def:J3_v2}
		\td J_3 := \la \na f(\bbx^{t,j_2-1}), \bbx^{t,j_2} - \bbv^{j_2} \ra 
	\end{equation}
	\begin{equation}\label{def:J4_v2}
		\td J_4 := \la \na f(\bbx^{t,M}) - \na f(\bbx^{t,j_1-1}), \bbv^{j_1} - \bbx^{t,j_2} \ra + 
		\la \na f(\bbx^{t,M}) - \na f(\bbx^{t,j_2-1}), \bbx^{t,j_2} - \bbv^{j_2} \ra
	\end{equation}
	
	We have the following lemma which gives upper bounds for $\td J_1, \td J_2, \td J_3$ and $ \td J_4 $.
	\begin{lemma}\label{lemma: J_1--J_5_v2}
		For any $t\ge 0$, $j_1 \in V^{t,M}$ and $j_2 \in [M]$ with $j_1\neq j_2$, the quantities defined in \eqref{def:J1-J2_v2}, \eqref{def:J3_v2} and \eqref{def:J4_v2} satisfy:
		
		(1) $\td J_1 \le LD \| \bbx^{t,j_1} - \bbx^{t,j_1-1} \| $. 
		
		(2) $\td J_2 \le \sum_{k=1}^M \la \na f(\bbx^{t,k-1}), \bbx^{t,k-1} - \bbx^{t,k} \ra + LD \sum_{k=1}^M \| \bbx^{t,k-1} - \bbx^{t,k} \| $. 
		
		(3) $\td J_3 \le  LD\| \bbx^{t,j_2-1} -  \bbx^{t,j_2} \|$.
		
		(4) $\td J_4 \le 2LD \sum_{k=1}^M\| \bbx^{t,k-1} - \bbx^{t,k} \| $. 
	\end{lemma}
	
	\noindent
	\textit{Proof of Lemma~\ref{lemma: J_1--J_5_v2}.}
		(1) We discuss different cases of $\al_{t,j_1} $. $(i)$
		If $\al_{t,j_1} = 1$, then $ \bbv^{j_1} = \bbx^{t,j_1}$ and hence $\td J_1 = 0$. 
		$(ii)$
		If $ \al_{t,j_1} \in
		(-\ga_{t,j_1} ,1)$, then by Lemma \ref{lemma: proximal-gradient-results-away} (c) we have $\td J_1 \le  LD \| \bbx^{t,j_1} - \bbx^{t,j_1-1} \|$. 
		$(iii)$ If $\al_{t,j_1} = -\ga_{t, j_1} $, 
		then $ \lam_{j_1}^{t,j_1} = 0 $, and by the updating rule \eqref{PolyCDwA update lam} we know $ \lam_{j_1}^{t,j_1+1} = \lam_{j_1}^{t,j_1+2} = \cdots = \lam_{j_1}^{t,M} = 0$. But this is contradictory to the fact that $ j_1 \in V^{t,M} = \{   j \in [M] ~|~ \lam^{t,M}_j >0    \}$, so the case $\al_{t,j_1} = -\ga_{t, j_1} $ cannot arise. 
		
		\smallskip
		\noindent
		(2) 
		From \eqref{shortcut-eq2} in Lemma~\ref{lemma: reduction} we know 
		\begin{equation}
			\begin{aligned}
				&\la \na f(\bbx^{t,j}) , \bbx^{t,j} - \bbz \ra -  	\la \na f(\bbx^{t,i-1}) , \bbx^{t,i} - \bbz \ra \\
				= \ &
				\sum_{k=i+1}^j   \la  \na f(\bbx^{t,k-1}), \bbx^{t,k} - \bbx^{t,k-1} \ra +   \sum_{k=i}^j 
				\la \na f(\bbx^{t,k})   - \na f(\bbx^{t,k-1})    ,   \bbx^{t,k} -\bbz \ra 	
			\end{aligned}
			\nonumber
		\end{equation}
		for all $i<j$ and $\bbz\in S$. This is equivalent to 
		\begin{equation}\label{shortcut-eq4}
			\begin{aligned}
				&\la \na f(\bbx^{t,j-1}) , \bbx^{t,j} - \bbz \ra -  	\la \na f(\bbx^{t,i-1}) , \bbx^{t,i} - \bbz \ra \\
				= \ &
				\sum_{k=i+1}^j   \la  \na f(\bbx^{t,k-1}), \bbx^{t,k} - \bbx^{t,k-1} \ra +   \sum_{k=i}^{j-1}
				\la \na f(\bbx^{t,k})   - \na f(\bbx^{t,k-1})    ,   \bbx^{t,k} -\bbz \ra 	
			\end{aligned}
		\end{equation}
		for all $i<j$ and $\bbz\in S$. 
		
		To bound $\td J_2$, we consider 2 cases. $(i)$ 
		If $j_2 <j_1$, then using \eqref{shortcut-eq4} with $i = j_2$ and $j = j_1$ and $\bbz = \bbx^{t,j_2}$, we have 
		\begin{equation}\label{shortcut-eq4pt5}
			\begin{aligned}
				\td J_2 =& 
				\la \na f(\bbx^{t,j_1-1}) , \bbx^{t,j_1} - \bbx^{t,j_2} \ra \\
				=&  \sum_{k=j_2+1}^{j_1}  \la  \na f(\bbx^{t,k-1}), \bbx^{t,k} - \bbx^{t,k-1} \ra  + 
				\sum_{k=j_2}^{j_1-1}
				\la \na f(\bbx^{t,k})   - \na f(\bbx^{t,k-1})    ,   \bbx^{t,k} -\bbx^{t,j_2} \ra. 	\\
				\le & \sum_{k=j_2}^{j_1-1}
				\la \na f(\bbx^{t,k})   - \na f(\bbx^{t,k-1})    ,   \bbx^{t,k} -\bbx^{t,j_2} \ra  \le \sum_{k=1}^{M}
				LD \| \bbx^{t,k}   - \bbx^{t,k-1}  \|
			\end{aligned}
			\nonumber
		\end{equation}
		where the first inequality is because $\la  \na f(\bbx^{t,k-1}), \bbx^{t,k} - \bbx^{t,k-1} \ra  \le0  $ for all $k\in [M]$ (by Lemma \ref{lemma: proximal-gradient-results-away} (a)).

		$(ii)$	If $j_1<j_2$, then using \eqref{shortcut-eq4} with $i = j_1$ and $j = j_2$ and $\bbz = \bbx^{t,j_2}$, we have 
		\begin{equation}
			\begin{aligned}
				\td J_2 = \ &
				\la \na f(\bbx^{t,j_1-1}) , \bbx^{t,j_1} - \bbx^{t,j_2} \ra \\
				= \ & \sum_{k=j_1+1}^{j_2}  \la  \na f(\bbx^{t,k-1}), \bbx^{t,k-1} - \bbx^{t,k} \ra  + 
				\sum_{k=j_1}^{j_2-1}
				\la \na f(\bbx^{t,k-1})   - \na f(\bbx^{t,k})    ,   \bbx^{t,k} -\bbx^{t,j_2} \ra \\
				\le \ & 
				\sum_{k=1}^M \la  \na f(\bbx^{t,k-1}), \bbx^{t,k-1} - \bbx^{t,k} \ra  + LD \sum_{k=1}^M  \| \bbx^{t,k-1}   - \bbx^{t,k}   \|	
			\end{aligned}
			\nonumber
		\end{equation}
		where the above inequality makes use of the fact $\la  \na f(\bbx^{t,k-1}), \bbx^{t,k-1} - \bbx^{t,k} \ra  \ge0  $ (by Lemma \ref{lemma: proximal-gradient-results-away} (a)). 
		Combining the cases $(i)$ and $(ii)$, we complete the proof of (2). 
		
		\smallskip
		\noindent
		(3) The conclusion can be immediately verified using Lemma~\ref{lemma: proximal-gradient-results-away} (b).

		\smallskip
		\noindent
		(4) Note that 
		\begin{equation}
			\begin{aligned}
				\td J_4  
				= \ & \la \na f(\bbx^{t,M}) - \na f(\bbx^{t,j_1-1}), \bbv^{j_1} - \bbx^{t,j_2} \ra + 
				\la \na f(\bbx^{t,M}) - \na f(\bbx^{t,j_2-1}), \bbx^{t,j_2} - \bbv^{j_2} \ra \\
				\le \ & D \| \na f(\bbx^{t,M}) - \na f(\bbx^{t,j_1-1})\| + D \|\na f(\bbx^{t,M}) - \na f(\bbx^{t,j_2-1}) \| \\
				\le \ &
				2D \sum_{k=1}^M\| \na f(\bbx^{t,k-1}) - \na f(\bbx^{t,k}) \| \le 
				2LD \sum_{k=1}^M\| \bbx^{t,k-1} - \bbx^{t,k} \|
			\end{aligned}
			\nonumber
		\end{equation}
		This completes of the proof of Lemma~\ref{lemma: J_1--J_5_v2}
\endproof

	The results of Lemma~\ref{lemma: bound J1-J5} immediately yields the following lemma. 
	\begin{lemma}\label{lemma:key1_v2}
		For any $t\ge 0$, $j_1 \in V^{t,M}$ and $j_2 \in [M]$ with $j_1\neq j_2$, it holds 
		\begin{equation}\label{h1}
			\la \na f(\bbx^{t,M}), \bbv^{j_1} - \bbv^{j_2} \ra
			\le
			\sum_{k=1}^M \la \na f(\bbx^{t,k-1}), \bbx^{t,k-1} - \bbx^{t,k} \ra + 4LD \sum_{k=1}^M \| \bbx^{t,k-1} - \bbx^{t,k} \| 
			\nonumber
		\end{equation}
		
	\end{lemma}
\textit{Proof of Lemma~\ref{lemma:key1_v2}.}
		By some algebra, 
		we have the decomposition:
		\begin{equation}
			\begin{aligned}
				&
				\la \na f(\bbx^{t,M}), \bbv^{j_1} - \bbv^{j_2} \ra \\
				= \ &
				\la \na f(\bbx^{t,M}), \bbv^{j_1} - \bbx^{t,j_1} \ra + 
				\la \na f(\bbx^{t,M}), \bbx^{t,j_1} - \bbx^{t,j_2} \ra +
				\la \na f(\bbx^{t,M}), \bbx^{t,j_2} - \bbv^{j_2} \ra  
				\\
				= \ &
				\td J_1 + \td J_2 + \td J_3 + \td J_4  .
			\end{aligned}
			\nonumber
		\end{equation}
		Using the above equality and Lemma~\ref{lemma: J_1--J_5_v2}, the conclusion is reached.
\endproof
	
	With Lemmas~\ref{lemma: proximal-gradient-results-away} and \ref{lemma:key1_v2} at hand, we are ready to present the proof of Theorem~\ref{theorem: PolyCDwA strong convex} (2). 
	Let $\wtd \mu := \mu \psi_S^2/4 $. Then 
	by Lemma \ref{lemma: key bound away}, there exist $\eta\in [0,2]$, $j_1 \in V^{t,M}$ and $j_2 \in [M]$ 
	such that 
	\begin{equation}\label{key3_v2}
		f(\bbx^{t,M}) - f(\bbx^*) 
		\le 
		\frac{\eta}{2}    \la \na f(\bbx^{t,M}), \bbv^{j_1} - \bbv^{j_2} \ra  - \frac{\wtd \mu}{2} \eta^2.
	\end{equation}
	Note that we can assume that $j_1 \neq j_2$, since otherwise $ f(\bbx^{t,M}) - f(\bbx^*) 
	\le 0$ and the conclusion of Theorem~\ref{theorem: PolyCDwA strong convex} (2) holds trivially. 
	First, by \eqref{key3_v2} and Lemma~\ref{lemma:key1_v2}, we have
	
	\begin{equation}\label{key5_v2}
		\begin{aligned}
			& f(\bbx^{t,M}) - f(\bbx^*)  \\
			\le \ &
			\frac{\eta}{2}    \sum_{k=1}^M \la \na f(\bbx^{t,k-1}), \bbx^{t,k-1} - \bbx^{t,k} \ra +
			2\eta LD \sum_{k=1}^M \| \bbx^{t,k} - \bbx^{t,k-1} \|   - \frac{\wtd \mu}{2} \eta^2 \\
			\le \ &
			\sum_{k=1}^M \la \na f(\bbx^{t,k-1}), \bbx^{t,k-1} - \bbx^{t,k} \ra +
			2\eta LD \sum_{k=1}^M \| \bbx^{t,k} - \bbx^{t,k-1} \|   - \frac{\wtd \mu}{2} \eta^2 
		\end{aligned}
	\end{equation}
	where the second inequality is because $\eta \le 2$. By Cauchy-Schwarz inequality we have 
	\begin{equation}\label{key6_v2}
		\begin{aligned}
			2\eta LD \sum_{k=1}^M \| \bbx^{t,k} - \bbx^{t,k-1} \| -  \frac{\wtd \mu}{2} \eta^2 
			\le \ &  (2\wtd \mu)^{-1}4L^2D^2 \Big(    \sum_{k=1}^M \| \bbx^{t,k} - \bbx^{t,k-1} \|  \Big)^2 \\
			\le \ &\frac{2ML^2D^2}{\wtd \mu}   \sum_{k=1}^M \| \bbx^{t,k} - \bbx^{t,k-1} \|^2
		\end{aligned}
	\end{equation}
	As a result, from \eqref{key5_v2} and \eqref{key6_v2} we have 
	\begin{equation}\label{h2}
		f(\bbx^{t,M}) - f(\bbx^*)  \le  \frac{2ML^2D^2}{\wtd \mu}   \sum_{k=1}^M \| \bbx^{t,k} - \bbx^{t,k-1} \|^2+  
		\sum_{k=1}^M \la \na f(\bbx^{t,k-1}), \bbx^{t,k-1} - \bbx^{t,k} \ra 
	\end{equation}
	On the other hand, by \eqref{ineq: convex-basic1} in Lemma~\ref{lemma: convex-basic} we have 
	\begin{equation}\label{h3}
		\begin{aligned}
			f(\bbx^{t,k-1}) - f(\bbx^{t,k}) \ge \ & \la \na f(\bbx^{t,k-1}), \bbx^{t,k-1} - \bbx^{t,k} \ra - \frac{L}{2} \|\bbx^{t,k-1} - \bbx^{t,k}  \|^2 \\
			\ge \ &
			(1/2)  \la \na f(\bbx^{t,k-1}), \bbx^{t,k-1} - \bbx^{t,k} \ra  \\
			\ge \ &
			(L/2) \|\bbx^{t,k-1} - \bbx^{t,k}  \|^2
		\end{aligned}
	\end{equation}
	where the second and third inequalities are both by Lemma \ref{lemma: proximal-gradient-results-away} (a). 
	Using the second inequality in~\eqref{h3} for $k\in [M]$, 
	\begin{equation}\label{h4}
		\sum_{k=1}^M  \la \na f(\bbx^{t,k-1}), \bbx^{t,k-1} - \bbx^{t,k} \ra  \le 2 ( f(\bbx^{t,0}) - f(\bbx^{t,M}) )
	\end{equation}
	and similarly, using the third inequality in~\eqref{h3}, we have:
	\begin{equation}\label{h5}
		\sum_{k=1}^M \|\bbx^{t,k-1} - \bbx^{t,k}  \|^2 \le \frac{2}{L} ( f(\bbx^{t,0}) - f(\bbx^{t,M}) ).
	\end{equation}
	
	Combining \eqref{h2}, \eqref{h4} and \eqref{h5}, we have 
	\begin{equation}
		f(\bbx^{t,M}) - f(\bbx^*)  \le  ({4ML D^2}/{\wtd \mu} + 2 )(f(\bbx^{t,0}) - f(\bbx^{t,M})) = 
		G' (f(\bbx^{t,0}) - f(\bbx^{t,M})). \nonumber
	\end{equation}
	Recall that $\bbx^{t,0} = \bbx^t$ and $\bbx^{t+1} = \bbx^{t,M}$, so we have 
	\begin{equation}
		f(\bbx^{t+1}) - f^* \le G' (f(\bbx^t) - f(\bbx^{t+1}))  \nonumber
	\end{equation}
	or equivalently $ f(\bbx^{t+1}) - f^* \le \frac{G'}{1+G'} ( f(\bbx^t) - f^* )  $, 
	which completes the proof.

	\section{Technical results}
	
	\begin{lemma}\label{lemma: convex-basic}
		[Theorem 2.1.5 of \cite{nesterov2003introductory}]
		Suppose $f(\cdot)$ is convex and $L$-smooth on $S$. Then for any $\bbx,\bby\in S$, it holds
		\begin{equation}\label{ineq: convex-basic1}
			f(\bby) \le f(\bbx)  + \la \na f(\bbx), \bby-\bbx \ra +  \frac{L}{2} \| \bby-\bbx\|^2
		\end{equation}
		and
		\begin{equation}\label{ineq: convex-basic2}
			f(\bby) \ge f(\bbx)  + \la \na f(\bbx), \bby-\bbx \ra +  \frac{1}{2L} \| \na f(\bby)- \na f(\bbx)\|^2
		\end{equation}
	\end{lemma}
	
	\begin{lemma}\label{lemma: technical1}
		Let $\{a_k\}_{k=1}^{\infty}$ be a sequence of decreasing positive numbers and $\lam >0$. Suppose $ a_k - a_{k+1} \ge \lam a_{k+1}^2$ for all $k\ge 1$. Then it holds 
		\begin{equation}
			a_k \le \frac{ \max\{a_1, 2/ \lam \} }{k} \quad \forall k\ge 1
			\nonumber
		\end{equation}
	\end{lemma}
\noindent
\textit{Proof of Lemma~\ref{lemma: technical1}.}
		Let $c :=\max\{a_1, 2/ \lam \}$. We prove the conclusion by induction. First, it is true that $ a_1 \le c = \frac{c}{1}$. Suppose the conclusion holds for $k$. From $ a_k - a_{k+1} \ge \lam a_{k+1}^2$ we have 
		\begin{eqnarray}\label{ineq1}
			a_{k+1} \le \frac{\sqrt{1+ 4\lam a_k} -1}{2\lam}          \le \frac{\sqrt{1+ 4\lam c/k } -1}{2\lam}     =      
			\frac{2c/k}{\sqrt{1+ 4\lam c/k} +1}
		\end{eqnarray}
		Since $c \ge  2/ \lam$, so we have 
		$ 2 \lam c -1 \ge \sqrt{1+4\lam c} \ge \sqrt{1+4\lam c/k}$
		for all $k\ge 1$. As a result, 
		\begin{equation}
			2 \le \frac{4\lam c}{ 1+ \sqrt{1+4\lam c /k}} = k ( \sqrt{1+4\lam c /k} - 1)
			\nonumber
		\end{equation}
		Hence $2 (k+1) \le k ( \sqrt{1+4\lam c /k} + 1)$. 
		Combining this inequality with \eqref{ineq1} we have 
		\begin{equation}
			a_{k+1} \le 	\frac{2c}{(\sqrt{1+ 4\lam c/k} +1)k} \le \frac{c}{k+1}
			\nonumber
		\end{equation}
		The proof is complete by induction. 
	\endproof

	\begin{lemma}\label{lemma: simplex vector decompose}
		For any $\bba,\bbb \in \Dt_M$, there exist $\bbp,\bbq\in \Dt_M$ such that 
		\begin{eqnarray}
			\bba-\bbb = \frac{\| \bba-\bbb\|_{1}}{2} (\bbp-\bbq) \quad {\rm{and}} \quad {\rm{supp}}(\bbp) \subseteq {\rm{supp}}(\bba) \ ,
			\nonumber
		\end{eqnarray}
		where ${\rm{supp}}(\bba)$ and ${\rm{supp}}(\bbp)$ denote the indices of nonzero coordinates of $\bba$ and $\bbp$ respectively. 
	\end{lemma}
\noindent
\textit{Proof of Lemma~\ref{lemma: simplex vector decompose}.}
		For any vector $\bbx\in \R^M$, let $\bbx^+$ be the vector in $\R^M$ with $x_i^+= \max \lt\{ x_i,0 \rt\}$ and $\bbx^- =  \bbx^+ - \bbx$.
		Assume $\bba\neq \bbb$ (otherwise the conclusion is trivial). 
		Let $p := 2 (\bba-\bbb)^+/\|\bba-\bbb\|_1$ and $q:= 2(\bba-\bbb)^-/\|\bba-\bbb\|_1$.
		Then it holds 
		$
		\bba-\bbb = \frac{\| \bba-\bbb\|_{1}}{2} (\bbp-\bbq)
		$~ and~ $\text{supp}(\bbp) \subseteq \text{supp}(\bba)$. Note that
		\begin{equation}
			\bbone_M^\T \bbp - \bbone_M^\T \bbq = \frac{2}{\|\bba-\bbb\|_1} (\bbone^\T_M \bba- \bbone_M^\T \bbb) = 0,~~~~~
			\bbone_M^\T \bbp + \bbone_M^\T \bbq = \frac{2}{\|\bba-\bbb\|_1} \|\bba-\bbb\|_1 = 2 \ .
			\nonumber
		\end{equation}
		As a result, it holds $\bbone_M^\T \bbp = \bbone_M^\T \bbq =1$, hence $\bbp,\bbq \in \Dt_M$. 
	\endproof

	\begin{lemma}\label{lemma: reduction}
		Let $\{ \bbx^{t,i}\}_{t \ge 0, 0\le  i  \le M} $ be any sequence in $S$. 
		Then for any 
		any $\bbz \in S$, for $0\le i<j \le M$ and $t\ge 0$, the following two equalities hold:
		\begin{equation}\label{shortcut-eq1}
			\begin{aligned}
				&\la \na f(\bbx^{t,j}) , \bbx^{t,j} -\bbz \ra -  	\la \na f(\bbx^{t,i}) , \bbx^{t,i} - \bbz \ra \\
				= \ &
				\sum_{k=i+1}^j \lt(  \la  \na f(\bbx^{t,k}), \bbx^{t,k} - \bbx^{t,k-1} \ra + \la \na f(\bbx^{t,k})   - \na f(\bbx^{t,k-1})    ,   \bbx^{t,k-1} -\bbz \ra \rt) 
			\end{aligned}
		\end{equation}
		and
		\begin{equation}\label{shortcut-eq2}
			\begin{aligned}
				&\la \na f(\bbx^{t,j}) , \bbx^{t,j} - \bbz \ra -  	\la \na f(\bbx^{t,i-1}) , \bbx^{t,i} - \bbz \ra \\
				= \ &
				\sum_{k=i+1}^j   \la  \na f(\bbx^{t,k-1}), \bbx^{t,k} - \bbx^{t,k-1} \ra +   \sum_{k=i}^j 
				\la \na f(\bbx^{t,k})   - \na f(\bbx^{t,k-1})    ,   \bbx^{t,k} -\bbz \ra 
			\end{aligned}
		\end{equation}
	\end{lemma}
	
	\noindent
	\textit{Proof of Lemma~\ref{lemma: reduction}.}
		Note that for any $i+1 \le k \le j$, 
		\begin{equation}
			\begin{aligned}
				&\la \na f(\bbx^{t,k}) , \bbx^{t,k} - \bbz \ra -  	\la \na f(\bbx^{t,k-1}) , \bbx^{t,k-1} - \bbz \ra \\
				= \ &
				\la \na f(\bbx^{t,k}) , \bbx^{t,k} - \bbz \ra - 	\la \na f(\bbx^{t,k}) , \bbx^{t,k-1} - \bbz \ra \\
				& + \la \na f(\bbx^{t,k}) , \bbx^{t,k-1} - \bbz \ra  - \la \na f(\bbx^{t,k-1}) , \bbx^{t,k-1} - \bbz \ra \\
				= \ &
				\la  \na f(\bbx^{t,k}), \bbx^{t,k}- \bbx^{t,k-1}  \ra + \la \na f(\bbx^{t,k})   - \na f(\bbx^{t,k-1})    ,   \bbx^{t,k-1} -\bbz \ra 	
			\end{aligned}
			\nonumber
		\end{equation}
		Summing it up from $i+1$ to $j$ we complete the proof of \eqref{shortcut-eq1}.  To prove \eqref{shortcut-eq2}, note that 
		\begin{equation}
			\begin{aligned}
				&\la  \na f(\bbx^{t,k}), \bbx^{t,k} - \bbx^{t,k-1} \ra + \la \na f(\bbx^{t,k})   - \na f(\bbx^{t,k-1})    ,   \bbx^{t,k-1} -\bbz \ra \\
				= \ &
				\la  \na f(\bbx^{t,k-1}), \bbx^{t,k} - \bbx^{t,k-1} \ra
				+ \la  \na f(\bbx^{t,k}) - \na f(\bbx^{t,k-1}), \bbx^{t,k} - \bbx^{t,k-1} \ra 
				\nonumber\\
				&+ \la \na f(\bbx^{t,k})   - \na f(\bbx^{t,k-1})    ,   \bbx^{t,k-1} -\bbz \ra \\
				= \ &
				\la  \na f(\bbx^{t,k-1}), \bbx^{t,k} - \bbx^{t,k-1} \ra + \la \na f(\bbx^{t,k})   - \na f(\bbx^{t,k-1})    ,   \bbx^{t,k} -\bbz \ra	
			\end{aligned}
			\nonumber
		\end{equation}
		for all $i+1 \le k \le j$. So by \eqref{shortcut-eq1} we know 
		\begin{equation}\label{shortcut-eq3}
			\begin{aligned}
				&\la \na f(\bbx^{t,j}) , \bbx^{t,j} - \bbz \ra -  	\la \na f(\bbx^{t,i}) , \bbx^{t,i} - \bbz \ra \\
				= \ &
				\sum_{k=i+1}^j \lt(  \la  \na f(\bbx^{t,k-1}), \bbx^{t,k} - \bbx^{t,k-1} \ra + \la \na f(\bbx^{t,k})   - \na f(\bbx^{t,k-1})    ,   \bbx^{t,k} -\bbz \ra \rt) 
			\end{aligned}
		\end{equation}
		Adding $\la \na f(\bbx^{t,i})   - \na f(\bbx^{t,i-1})    ,   \bbx^{t,i} -\bbz \ra$ in both sides of \eqref{shortcut-eq3} we complete the proof of \eqref{shortcut-eq2}. 
	\endproof

	\begin{lemma}\label{lemma: poly}
		Let $S$ be a bounded polyhedron and $\psi_S$ be defined in \eqref{facial-distance-defn}. 
		For any $\bbx,\bby\in S$, given $\bblam^x\in \Dt_M$ satisfying $ \bbx = \sum_{i=1}^M \lam^x_{i} \bbv^i  $, there exists $\bblam^y \in \Dt_M$ such that $\bby = \sum_{i=1}^M \lam^y_{i} \bbv^i  $ and $\|\bbx-\bby\| \ge ({\psi_S}/{2})\| \bblam^x - \bblam^y \|_1 $.
	\end{lemma}
	The proof of Lemma \ref{lemma: poly} can be found in the proof of \cite{gutman2018condition} Proposition 1.

	\section{Additional experimental details}

	\subsection{Comparison of PolyCDwA with benchmarks}\label{app-subsubsection: linear regression Comparison with benchmarks}
	Implementation details of benchmark methods:

	\begin{itemize}
		
		\item {\bf 2-CD}: Randomized 2-coordinate methods following~\cite{reddi2014large,necoara2017random}. This algorithm is applied to the reformulation of \eqref{problem: least squares l1 constraints} as an optimization problem on the standard simplex.
		In each iteration, we select (uniformly at random) two coordinates to swap. The algorithm is run with a maximum of $100d$ iterations (where $d$ is the dimension of $x$).
		
		\item {\bf AFW}: Away-step FW method \cite{lacoste2015global}. This method is run for a maximum of $5000$ iterations  and is terminated earlier if the relative improvement in the past $50$ iterations is less than $10^{-8}$, that is, if we let $f_k$ be the best objective value found in the first $k$ iterations, 
		then it is terminated at iteration $k$ if
		$(f_{k-50} - f_k)/\max\{|f_{k-50}|,1\}<10^{-8}$. 
		
		\item {\bf FISTA}: The accelerated proximal gradient method in \cite{beck2009fast}. We run FISTA for a maximum of $1000$ iterations.
		
		\item {\bf StrOpt}: This is used with the default settings of the Julia package StructuredOptimization.jl. 
		
		\item {\bf MOSEK}: The interior-point commercial solver MOSEK \cite{andersen2000mosek}. This method is called with the Julia interface Convex.jl and used with the default setting. 
	\end{itemize}

	\subsection{Kernel density estimation}\label{app-subsection: exp: kernal density estimation}
	Data generation:
	Given ${\mu}\in \R^d$ and $\Sigma \in \mathbb{S}^{d\times d}_+$, 
	let $\eta(\bbmu, \Sigma, \cdot)$ denote the density function of normal distribution $N(\bbmu, \Sigma)$. 
	Let $g^*(\cdot)$ be a Gaussian mixture density with $m =10$ component $g^*(\bbx) = \sum_{j=1}^{m} \bar w_j \cdot \eta(\bar\bbmu_j, \bar\sigma_j I_d , x) $, where $\bar \bbw = [\bar w_1,...,\bar w_m]^\T $ is generated uniformly randomly on $ \Dt_m$; 
	$\bar\bbmu_j$ ($j\in [m]$) are iid from $N({0}_d,16I_d)$;
	$\bar\sigma_j$ ($j\in [m]$) are iid uniformly randomly generated from $[0.2,1.2]$. 
	Given $n\ge1 $ and denote $n_0 = \lfloor n/100 \rfloor$ and $n_1 = n - n_0$. Then we generate $n$ samples $X_1,...., X_n$, where $n_1$ of them are iid from $g^*(\cdot)$, and the remaining $n_0$ are outliers generated iid from $N({0}_d, 2500 I_d)$. 
	For the experiment in Figure~\ref{fig: density compare}, we take
	$\sigma = 1.0$ (in \eqref{gaussian-kernal}) and $\mu = 0.4$ (in \eqref{Huber-loss}).

	\section*{Acknowledgments} The authors would like to thank Robert Freund and Zikai Xiong for their helpful discussions and comments on the paper.

	\nocite{lacoste2013block}
	\nocite{wang2016parallel}
	\nocite{beck2015block}
	\nocite{CGFWSurvey2022}

	\bibliographystyle{plain}     
	\bibliography{mybib}

\end{document}